\documentclass[12pt]{article}
\usepackage{amsmath,amssymb,amsfonts,fancybox,ascmac,dashbox,color}
\usepackage{graphicx}
\usepackage{mathtools}
\mathtoolsset{showonlyrefs=true}
\usepackage{tikz,braket}
\usetikzlibrary{intersections,calc}

\newcommand{\R}{\mathbb{R}}
\newcommand{\N}{\mathbb{N}}

\newcommand{\I}{\mathcal{I}}
\newcommand{\T}{\mathcal{T}}
\newcommand{\LL}{\mathcal{L}}
\newcommand{\PP}{\mathcal{P}}
\newcommand{\dd}{\text{d}}
\newcommand{\lam}{\lambda}
\newcommand{\bfx}{\mathbf{x}}
\newcommand{\bfy}{\mathbf{y}}
\newcommand{\bfw}{\mathbf{w}}

\newcommand{\bfb}{\mathbf{b}}
\newcommand{\hK}{\widehat{K}}

\newcommand{\Fab}{F_{\alpha\beta}}
\newcommand{\Kab}{K_{\alpha\beta}}
\newcommand{\wA}{\widetilde{A}}
\newcommand{\hv}{\hat{v}}
\newcommand{\Simp}{\mathbb S}
\newcommand{\vvskip}{\vspace{0.3truecm}}

\newtheorem{theorem}{Theorem}
\newtheorem{lemma}[theorem]{Lemma}
\newtheorem{corollary}[theorem]{Corollary}

\newtheorem{assumption}[theorem]{Assumption}

\begin{document}
%
%
\title{Lectures on the Error Analysis of Interpolation \\ 
on Simplicial Triangulations without \\[5pt]
the Shape Regularity Assumption \\[5pt]
Part 1: Lagrange Interpolation on Triangles
}
\author{Kenta Kobayashi
\footnote{Graduate School of Business Administration, 
         Hitotsubashi University, Kunitachi, JAPAN} \quad
Takuya Tsuchiya
\footnote{Graduate School of Science and Engineering, Ehime University,
Matsuyama, JAPAN, \newline
\hspace*{6mm} \texttt{tsuchiya@math.sci.ehime-u.ac.jp}.}}

\date{January 18, 2022}
\setlength{\baselineskip}{18pt}

\maketitle

\noindent \textbf{Abstract:}
In the error analysis of finite element methods, the shape regularity
assumption on triangulations is typically imposed to obtain 
\textit{a priori} error estimations.  In practical computations,
however, very ``thin'' or ``degenerated'' elements that violate the
shape regularity assumption may appear when we use adaptive mesh
refinement.  In this survey, we attempt to establish an error
analysis approach without the shape regularity assumption on
triangulations.

We have presented several papers on the error analysis of
finite element methods on non-shape regular triangulations.
The main points in these papers are that,
\textit{in the error estimates of finite element methods, the
circumradius of the triangles is one of the most important factors}.

The purpose of this survey is to provide a simple and plain
explanation of the results to researchers and, in particular,
graduate students who are interested in the subject.
Therefore, this survey is not intended to be a research paper.
We hope that, in the near future, it will be merged into a textbook on
the mathematical theory of the finite element methods.

\section{Introduction: Lagrange interpolation on triangles}

Lagrange interpolation on triangles and the associated error
estimates are important subjects in numerical analysis. In particular,
they are crucial in the error analysis of finite element methods.
Throughout this survey, $K \subset \R^2$ denotes a triangle
with vertices $\bfx_i$, $i=1,2,3$.  In this survey, we always assume
that triangles are closed sets.  Let $\lam_i$ be the barycentric
coordinates of $K$ with respect to $\bfx_i$. By definition,
$0 \le \lam_i \le 1$, $\sum_{i=1}^{3} \lam_i =1$.  Let $\N_0$
be the set of nonnegative integers, and
$\gamma = (a_1,a_2,a_3) \in \N_0^{3}$ be a multi-index.
Let $k$ be a positive integer. If $|\gamma| := \sum_{i=1}^{d+1}a_i = k$, then
${\gamma}/{k} := \left(a_1/k,a_2/k,a_3/k\right)$
can be regarded as a barycentric coordinate in $K$.
The set $\Sigma^k(K)$ of points on $K$ is defined as
\footnote{The set $\Sigma^k(K)$ is sometimes called a \textit{stencil}.}
\begin{equation}
   \Sigma^k(K) := \left\{
   \frac{\gamma}{k} \in K \Bigm| |\gamma| = k, \;
     \gamma \in \N_0^{3} \right\}.
   \label{Sigma}
\end{equation}

\begin{figure}[thb]
\begin{center}
  \includegraphics[width=10truecm]{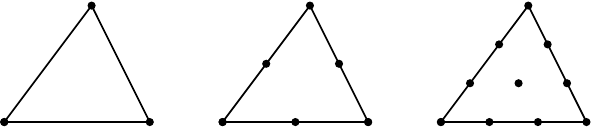} 
\caption{Set $\Sigma^k(K)$, $k=1$, $k=2$, $k=3$.}
 \label{gridpoints}
\end{center}
\end{figure}

Let $\PP_k(K)$ be a set of polynomials defined on $K$
whose degree is at most $k$.  For a continuous function
$v \in C^0(K)$, the $k$th-order Lagrange interpolation
$\I_K^k v \in \PP_k(K)$ is defined as
\begin{align*}
   v(\bfx) = (\I_K^kv)(\bfx), \qquad \forall \bfx \in \Sigma^k(K).
\end{align*}

To enable the error analysis of Lagrange interpolation, we typically
introduce the following condition
 \cite{Ciarlet, BrennerScott, ErnGuermond}.
Let $h_K := \mathrm{diam}K$ and $\rho_K$ be the diameter of its
inscribed circle.  Suppose that $X$ is a set of (possibly infinitely
many) triangles.
\begin{center}
\fbox{
\begin{minipage}{15truecm}
\begin{assumption}[Shape regularity]
The set $X$ is called \textbf{shape regular} if there exists
a constant $\sigma > 0$ such that
\begin{align*}
  \frac{h_K}{\rho_K} \le \sigma, \qquad \forall K \in X.
\end{align*}

\end{assumption}
\end{minipage}
}
\end{center}
The maximum of the ratio $h_K/\rho_K$ in $X$ is called its
\textbf{chunkiness parameter} \cite{BrennerScott}.  The shape regularity
condition is sometimes also called the \textbf{inscribed ball
condition}. For more information on the conditions equivalent to
shape regularity, see \cite{BrandtsKorotovKrizek}.

Let $\hK$ be a \textbf{reference element}.  The triangle with vertices
$(0,0)^\top$, $(1,0)^\top$, and $(0,1)^\top$ is typically taken as the
reference triangle $\hK$.   Let
$\varphi(\bfx) = A\bfx + \bfb$ be an affine transformation that maps
$\hK$ to $K$, where $A$ is a $2 \times 2$ regular matrix and
$\bfb \in \R^2$. 

Error analysis is first performed on the reference element
$\hK$.  Then, the ``pull back'' with $v \circ \varphi$ is used
to transfer the result obtained on $\hK$ to the ``physical
element'' $K$.

Let $\|A\|$ denote the matrix norm of $A$ associated with the Euclidean
norm of $\R^2$, and let $1 \le p \le \infty$.  The function
$v \in W^{k+1,p}(K)$ is pulled back by $\varphi$ as
$\hat{v} := v \circ \varphi$.  Let $k$ and $m$ be integers such that
$k \ge 1$ and $0 \le  m \le k$.  The following theorem is standard.
\begin{center}
\fbox{
\begin{minipage}{15truecm}
\begin{theorem}[\cite{Ciarlet}, Theorem~3.1.4]\label{Thm314}
Let $\sigma > 0$ be a constant. If $h_K/\rho_K \le \sigma$, then there
exists a constant $C = C(\hK,p,k,m)$ independent of $K$ such that,
for $v \in W^{k+1,p}(K)$,
\begin{align}
  |v - \I_K^k v|_{m,p,K} & \le C \|A\|^{k+1}\|A^{-1}\|^m
   |v|_{k+1,p,K} \notag \\
  & \le C \frac{h_K^{k+1}}{\rho_K^m} |v|_{k+1,p,K}
   \le (C\sigma^m) h_K^{k+1-m} |v|_{k+1,p,K}.
  \label{standard-est0}
\end{align}
\end{theorem}
\end{minipage}
}
\end{center}

To derive the second inequality in \eqref{standard-est0}, we use
the following lemma.
%
\begin{center}
\fbox{
\begin{minipage}{15truecm}
\begin{lemma}[\cite{Ciarlet}, Theorem~3.1.3]\label{inscribed}
 We have $\|A\| \le h_K \rho_{\hK}^{-1}$, 
$\|A^{-1}\| \le h_{\hK} \rho_K^{-1}$.
\end{lemma}
\end{minipage}
}
\end{center}

Let $K$ be an arbitrary triangle, and $h_K \ge \alpha \ge \beta > 0$
be the lengths of its three edges.  Note that $h_K/2 < \alpha \le h_K$.
Using translation, rotation, and mirror imaging, $K$ is transformed into 
a triangle with vertices 
$\bfx_1 = (0,0)^\top$, $\bfx_2 = (\alpha,0)^\top$, and
$\bfx_3 = (\beta s, \beta t)^\top$, where $s = \cos\theta$,
$t = \sin\theta$, and $0 < \theta < \pi$ is the inner angle of $K$ at
$\bfx_1$.  This triangle is called the \textbf{standard position} of $K$.
By the law of cosines, 
\begin{align*}
   h_K^2 = \alpha^2 + \beta^2 - 2 \alpha \beta \cos \theta \quad
 \text{ and } \quad
  \cos \theta = \frac{\beta}{2\alpha} +
  \frac{\alpha^2 - h_K^2}{2\alpha\beta} 
  \le \frac{\beta}{2\alpha} \le \frac{1}{2}.
\end{align*}
Hence, $\pi/3 \le \theta < \pi$.

\begin{figure}[thb]
\begin{center}
\begin{tikzpicture}[line width = 1pt,scale=0.6]
   \coordinate [label=below:{$\bfx_1$}](A) at (0.0,0.0);
   \coordinate [label=below:{$\bfx_2$}](B) at (7.0,0.0);
   \coordinate [label=above:{$\bfx_3$}](C) at (-2.0,4.0);
   \draw (A) -- node[below]{$\alpha$} (B) ;
   \draw (B) -- node[pos=0.4,above]{$h_K$} (C) ;
   \draw (C) -- node[pos=0.6, left]{$\beta$} (A) ;
   \coordinate (D) at (0.8,0.0);
   \coordinate (E) at ($(A)!0.17!(C)$);
   \draw [bend right,thin] (D) to node[pos=0.3,above]{$\theta$} (E) ;
   \coordinate [label=:{$K$}](F) at (2.0,0.57);
\end{tikzpicture}
 \caption{General triangle $K$ in the standard position. 
 The vertices are
 $\bfx_1=(0,0)^\top$, $\bfx_2=(\alpha,0)^\top$, and
  $\bfx_3=(\beta s,\beta t)^\top$,
 where $s^2 + t^2 = 1$, $t > 0$.
 We assume that $0 < \beta \le \alpha \le h_K$.}
 \label{LiuKikuchiTriangle0}
\end{center}
\end{figure}
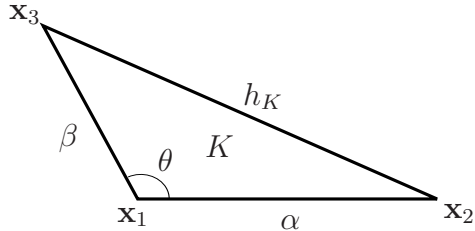
These assumptions imply that the affine transformation $\varphi$
can be written as $\varphi(\bfx) = A\bfx$ with the matrix
\begin{align}
  A = \begin{pmatrix}
      \alpha & \beta s \\
       0     & \beta t
      \end{pmatrix}.
   \label{mat-A}
\end{align}
We set $t = \sin\theta = 1$, for example (i.e., $K$ is a right triangle).
Then, $s = 0$, $\|A\| = \alpha$, $\|A^{-1}\| = 1/\beta$, and the inequalities
in \eqref{standard-est} can be rearranged as
\begin{align}
  |v - \I_K^k v|_{m,p,K} \le C \frac{\alpha^{k+1}}{\beta^m}
   |v|_{k+1,p,K}  \le C
  \left(\frac{\alpha}{\beta}\right)^m h_K^{k+1-m} |v|_{k+1,p,K}.
   \label{ab-est}
\end{align}
Thus, we might consider that the ratio $\alpha/\beta$ should not be
too large, or $K$ should not be too ``flat.''  This consideration is
expressed as the \textit{minimum angle condition}
(Zl\'amal \cite{Zlamal}, \v{Z}en\'i\v{s}ek \cite{Zenisek}), which
is equivalent to the shape regularity condition for triangles.

\vvskip
\begin{center}
\fbox{
\begin{minipage}{15truecm}
\begin{theorem}[Minimum angle condition]
Let $\theta_0$, $(0 < \theta_0 \le \pi/3)$ be a constant. If
any angle $\theta$ of $K$ satisfies $\theta \ge \theta_0$
and $h_K \le 1$, then there exists a constant $C = C(\theta_0)$
independent of $h_K$ such that
\begin{align*}
   |v - \I_K^1 v|_{1,2,K} \le C h_K |v|_{2,2,K}, \qquad
   \forall v \in H^2(K).
\end{align*}
\end{theorem}
\end{minipage}
}
\end{center}

However, the minimum angle condition and shape regularity are not
necessarily needed to obtain an error estimate.  The following condition
is well known (Babu\v{s}ka--Aziz \cite{BabuskaAziz}).
\begin{center}
\fbox{
\begin{minipage}{15truecm}
\begin{theorem}[Maximum angle condition]
Let $\theta_1$, $(\pi/3 \le \theta_1 < \pi)$ be a constant.  If
any angle $\theta$ of $K$ satisfies $\theta \le \theta_1$ and
$h_K \le 1$, then there exists a constant $C = C(\theta_1)$ 
that is independent of $h_K$ such that
\begin{equation}
   |v - \I_K^1 v|_{1,2,K} \le C h_K |v|_{2,2,K}, \qquad
   \forall v \in H^2(K).  
 \label{maximumangle}
\end{equation}
\end{theorem}
\end{minipage}
}
\end{center}

\vvskip
K\v{r}\'{i}\v{z}ek \cite{Krizek1} introduced
the \textit{semiregularity condition}, which is equivalent to
the maximum angle condition (see \textit{Remark} below).
Let $R_K$ be the circumradius of $K$.
\begin{center}
\fbox{
\begin{minipage}{15truecm}
\begin{theorem}[Semiregularity condition]
Let $p > 1$ and  $\sigma > 0$ be a constant. If $R_K/h_K \le \sigma$ and
$h_K \le 1$, then there exists a constant $C = C(\sigma)$ 
that is independent of $h_K$ such that
\begin{align*}
   |v - \I_K^1 v|_{1,p,K} \le C h_K |v|_{2,p,K}, \qquad
   \forall v \in W^{2,p}(K).
\end{align*}
\end{theorem}
\end{minipage}
}
\end{center}

\vvskip
We mention a few more known results.  Jamet \cite{Jamet} presented 
the following results.

\begin{center}
\fbox{
\begin{minipage}{15truecm}
\begin{theorem}\label{Jamettheorem}
Let $1 \le p \le \infty$. Let $m \ge 0$, $k \ge 1$ be integers such that
$k+1-m > 2/p$ $(1 < p \le \infty)$ or $k-m \ge 1$ $(p=1)$.
Then, the following estimate holds:
\begin{equation}
  |v - \I_K^k v|_{m,p,K} \le C \frac{h_K^{k+1-m}}{\cos^m \theta_K/2}
  |v|_{k+1,p,K},   \quad \forall v \in W^{k+1,p}(K),
  \label{jamet}
\end{equation}
where $\theta_K$ is the maximum angle of $K$, and $C$ depends only on
$k$ and $p$.
\end{theorem}
\end{minipage}
}
\end{center}

\noindent
\textsl{Remark:}
(1) In Theorem~\ref{Jamettheorem}, the restriction on $p$ comes from the
Sobolev imbedding theorem.  Note that in \cite[Th\'{e}or\`{e}me~3.1]{Jamet}
the case $p=1$ is not mentioned explicitly but clearly holds for
triangles (see Section~\ref{Sobolevimbedding}). 
For the case of the maximum angle condition, we set $k=m=1$ and find
that Jamet's result (Theorem~\ref{Jamettheorem}) \textit{does not} imply
the estimation \eqref{maximumangle} because the case $p = 2$ is excluded.

\noindent
(2) Let an arbitrary triangle $K$ be in its standard position
(Figure~\ref{LiuKikuchiTriangle0}).  Then $\theta$ is the maximum
internal angle of $K$, and
\begin{align}
    \frac{R_K}{h_K}  = \frac{1}{2\sin\theta}, \qquad
     \frac{\pi}{3} \le \theta < \pi 
     \label{semi-regular}
\end{align}
by the law of sines.  Thus, the dimensionless quantity 
$R_K/h_K$ represents the maximum internal angle of $K$, and
the boundedness of $R_K/h_K$, which
is the semiregularity of $K$, 
is equivalent to
the maximum angle condition $\theta \le \theta_1 < \pi$ with a
fixed constant $\theta_1$.  $\square$

\vvskip
For further results of the error estimations on ``skinny elements'',
see the monograph by Apel \cite{Apel}.

\vvskip
Recently, Kobayashi, one of the authors, obtained the following
epoch-making result \cite{Kobayashi1}. Let $A$, $B$, and $C$ be the
lengths of the three edges of $K$ and $S$ be the area of $K$.

\begin{center}
\fbox{
\begin{minipage}{15truecm}
\begin{theorem}[Kobayashi's formula]
We define the constant $C(K)$ as
\begin{align*}
   C(K) := \sqrt{\frac{A^2B^2C^2}{16S^2} - 
                    \frac{A^2 + B^2 + C^2}{30} - \frac{S^2}5
       \left(\frac1{A^2} + \frac1{B^2} + \frac1{C^2}\right)}.
\end{align*}
Then the following holds:
\begin{align*}
   |v - \I_K^1 v|_{1,2,K} \le C(K) |v|_{2,2,K}, \qquad
   \forall v \in H^2(K).
\end{align*}
\end{theorem}
\end{minipage}
}
\end{center}

\vvskip
Recall that $R_K$ is the circumradius of $K$ and is written as
\footnote{This formula is proved using the law of sines.}
\begin{align}
    R_K = \frac{ABC}{4S}.
   \label{circumradius-f}
\end{align}
Then, we immediately realize that
$C(K) < R_K$ and obtain a corollary of Kobayashi's formula.
\vvskip
\begin{center}
\fbox{
\begin{minipage}{15truecm}
\begin{corollary}
For any triangle $K \subset \R^2$, the following estimate holds:
\begin{equation}
    |v - \I_K^1 v|_{1,2,K} \le R_K |v|_{2,2,K}, \qquad
   \forall v \in H^2(K). 
   \label{Kobayashi-corollary}
\end{equation}
\end{corollary}
\end{minipage}
}
\end{center}

\begin{figure}[thb]
\begin{center}
  \includegraphics[width=5cm]{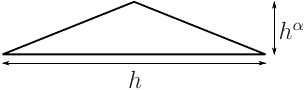} \qquad\quad
  \includegraphics[width=4.8cm]{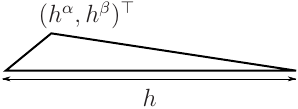}  \\[0.0cm]
\caption{Examples of triangles that violate the maximum
angle condition but satisfy $R_K \to 0$ as $h\to 0$.}
 \label{example1}
\end{center}
\end{figure}

\vvskip
This corollary demonstrates that even if the minimum angle is very small
or the maximum angle is very close to $\pi$, the error
$|v - \I_K^1 v|_{1,K}$ converges to $0$ if $R_K$ converges to $0$.
We consider the isosceles triangle $K$ shown in
Figure~\ref{example1} (left).   Using \eqref{circumradius-f}, we realize
that $R_K = h^\alpha/2 + h^{2-\alpha}/8 = \mathcal{O}(h^{2-\alpha})$
($\alpha \ge 1$, $h \le 1$). Thus, if $\alpha < 2$, $R_K \to 0$ as
$h \to 0$.

As another example, let $\alpha$, $\beta \in \R$ satisfy
$1 < \alpha < \beta < 1 + \alpha$.  
We consider the triangle $K$ whose vertices are
$(0,0)^\top$, $(h,0)^\top$, and $(h^\alpha,h^\beta)^\top$
(Figure~\ref{example1} (right)).
With \eqref{circumradius-f}, it is straightforward to see
\begin{align}
  R_K & = \frac{h
    \left(h^{2\alpha} + h^{2\beta}\right)^{1/2}
    \left((h^{\alpha} - h)^2 + h^{2\beta}\right)^{1/2}}
    {2h^{1+\beta}} \\
   & = \frac{h^{1 + \alpha}}{2h^{\beta}}
    \left(1 + h^{2\beta - 2\alpha}\right)^{1/2}
    \left(1 + h^{2\alpha - 2} - 2 h^{\alpha -1} 
    + h^{2\beta-2}\right)^{1/2} = \mathcal{O}(h^{1+\alpha - \beta}), \\
  \rho_K & = \frac{h^{1+\beta}}{h +
    \left(h^{2\alpha} + h^{2\beta}\right)^{1/2}
   + \left((h^{\alpha} - h)^2 + h^{2\beta}\right)^{1/2}}, \quad
   \frac{h^{\beta}}{3} < \rho_K < h^{\beta}.
\end{align}
Hence, if $h \to 0$, the convergence rates that
\eqref{standard-est0} and \eqref{Kobayashi-corollary} yield are
$\mathcal{O}(h^{2-\beta})$ and $\mathcal{O}(h^{1+\alpha - \beta})$,
respectively.  Therefore, \eqref{Kobayashi-corollary}
obtains a better convergence rate than \eqref{standard-est0}.  Moreover,
if $\beta \ge 2$, \eqref{standard-est0} does not yield convergence
whereas \eqref{Kobayashi-corollary} does.
Note that, when $h \to 0$, the maximum angles of $K$ approach to $\pi$
in both cases.

\vvskip

Although Kobayashi's formula is remarkable, its proof is long
and needs validated numerical computation.  We began this 
research to provide a ``paper-and-pencil'' proof of
\eqref{Kobayashi-corollary}, and recently reported an error
estimation in terms of the circumradius of a triangle
\cite{KobayashiTsuchiya1, KobayashiTsuchiya3, KobayashiTsuchiya4}.

\begin{center}
\fbox{
\begin{minipage}{15truecm}
\begin{theorem}[Circumradius estimates]\label{thm:circumradius-est}
Let $K$ be an arbitrary triangle.  Then, for the $k$th-order Lagrange
interpolation $\I_K^k$ on $K$, the estimation
\begin{equation} \label{eq:circumradius-est}
   |v - \I_K^k v|_{m,p,K} \le C
   \left(\frac{R_K}{h_K}\right)^m h_K^{k+1-m} |v|_{k+1,p,K}
  =  C R_K^m h_K^{k+1-2m} |v|_{k+1,p,K}
\end{equation}
holds for any $v \in W^{k+1,p}(K)$,
where the constant $C=C(k,m,p)$ is independent of the geometry of $K$.
\end{theorem}
\end{minipage}
}
\end{center}

We recall that a general triangle $K$ may be written using the settings
in Figure~\ref{LiuKikuchiTriangle0}.  The essence of the proof of 
Theorem~\ref{thm:circumradius-est} is that the matrix $A$ in
\eqref{mat-A} is decomposed as
\begin{align*}
   A = \wA D_{\alpha\beta}, \qquad
   \wA := \begin{pmatrix} 1 & s \\ 0 & t \end{pmatrix},
   \quad
   D_{\alpha\beta} := \begin{pmatrix} \alpha & 0 \\ 0 & \beta
       \end{pmatrix}.
\end{align*}
With this decomposition, the estimate \eqref{standard-est0} is
rearranged as
\begin{align}
  |v - \I_K^k v|_{m,p,K} & \le C \|\widetilde{A}\|^{k+1}
   \|\widetilde{A}^{-1}\|^{m}
   \|D_{\alpha\beta}\|^{k+1}\|D_{\alpha\beta}^{-1}\|^m
   |v|_{k+1,p,K}.
  \label{standard-est1}
\end{align}
As indicated by us \cite{KobayashiTsuchiya4} and Babu\v{s}ka--Aziz
\cite{BabuskaAziz}, the linear transformation by
$D_{\alpha\beta}$ does not reduce the approximation property of Lagrange
interpolation, and only $\wA$ could make it ``bad.''  This means
that the term
\begin{align*}
  \|D_{\alpha\beta}\|^{k+1}\|D_{\alpha\beta}^{-1}\|^m
    = \frac{(\max\{\alpha,\beta\})^{k+1}}
       {(\min\{\alpha,\beta\})^{m}}  \quad
\text{may be replaced with} \; C_1 h_K^{k+1-m}.
\end{align*}
Furthermore, $\|\wA\|$ and $\|\wA^{-1}\|$ (the maximum singular values
of $\wA$ and $\wA^{-1}$) are bounded using the circumradius $R_K$
and $h_K$ as
\begin{align*}
  \|\widetilde{A}\|^{k+1}\|\widetilde{A}^{-1}\|^m
   \le C_2 \left(\frac{R_K}{h_K}\right)^m, \quad
   \frac{R_K}{h_K} = \frac{1}{2\sin \theta},
\end{align*}
where $\theta$ is the maximum internal angle of $K$
(see Figure~\ref{LiuKikuchiTriangle0} and \eqref{semi-regular}).
 We emphasize that the constants $C_i$ $(i=1,2)$ only depend
on $k$, $m$, and $p$.  Note that, by setting $t = 1$ and
$\beta = \alpha^2$ in \eqref{standard-est0} (and \eqref{ab-est}), we
realize that, regardless of how much we try to analyze
$\|A\|^{k+1}\|A^{-1}\|^m$, we cannot prove 
Theorem~\ref{thm:circumradius-est}.  
In the sequel of this survey, we will explain the proof of
Theorem~\ref{thm:circumradius-est} in detail.

\section{Preliminaries}
\subsection{Notation}
Let $n \ge 1$ be a positive integer and $\R^n$ be $n$-dimensional
Euclidean space.  We denote the Euclidean norm of $\bfx \in \R^n$ by
$|\bfx|$.  Let $\R^{n*} := \{l:\R^n \to \R : l \text{ is linear}\}$
be the dual space of $\R^n$.  We always regard $\bfx \in \R^n$ as a
column vector and $\mathbf{a} \in \R^{n*}$ as a row vector.
For a matrix $A$ and $\bfx \in \R^n$, $A^\top$ and $\bfx^\top$
denote their transpositions.
For matrices $A = (a_{ij})_{i,j=1, \cdots, n}$ and
$B = (b_{ij})_{i,j=1,\cdots,n}$, their Kronecker product 
$A \otimes B$ is an $n^2 \times n^2$ matrix defined as
\begin{align*}
 A \otimes B := \begin{pmatrix}
      a_{11}B & \cdots & a_{1n}B \\
        \vdots & & \vdots \\
      a_{n1}B & \cdots & a_{nn}B \\
     \end{pmatrix}.
\end{align*}
For matrices $A_i$, $i = 1, \cdots, k$, the Kronecker product
$A_1 \otimes \cdots \otimes A_k$ is defined recursively.

For a differentiable function $f$ with $n$ variables,
its gradient $\nabla f = \mathrm{grad} f \in \R^{n*}$ is the row vector
defined as
\[
  \nabla f = \nabla_\bfx f := 
  \left(\frac{\partial f}{\partial x_1}, \cdots , 
    \frac{\partial f}{\partial x_n}\right), \qquad
    \bfx := (x_1, \cdots , x_n)^\top.
\]

Let $\N_{0}$ be the set of nonnegative integers.
For $\delta = (\delta_1,...,\delta_n) \in (\N_{0})^n$,
the multi-index $\partial^\delta$ of partial differentiation 
(in the sense of distribution) is defined by
\[
    \partial^\delta = \partial_\bfx^\delta
    := \frac{\partial^{|\delta|}\ }
   {\partial x_1^{\delta_1} \cdots \partial x_n^{\delta_n}}, \qquad
   |\delta| := \delta_1 + \cdots + \delta_n.
\]
For two multi-indices  $\eta=(\eta_1,\cdots,\eta_n)$,
$\delta=(\delta_1,\cdots,\delta_n)$,
$\eta \le \delta$ means that $\eta_i \le \delta_i$ $(i=1,\cdots,n)$.
Additionally, $\delta\cdot\eta$ and $\delta !$ are defined as
$\delta\cdot\eta := \eta_1\delta_1 + \cdots + \eta_n\delta_n$ and
$\delta ! := \delta_1! \cdots \delta_n!$, respectively.

Let $\Omega \subset \R^n$ be a (bounded) domain.  The usual Lebesgue
space is denoted by $L^p(\Omega)$ for $1 \le p \le \infty$.
For a positive integer $k$, the Sobolev space $W^{k,p}(\Omega)$ is
defined by
$\displaystyle
  W^{k,p}(\Omega) := 
  \left\{v \in L^p(\Omega) \, | \, \partial^\delta v \in L^p(\Omega), \,
   |\delta| \le k\right\}$.
For $1 \le p < \infty$, the norm and semi-norm of $W^{k,p}(\Omega)$ are
defined as
\begin{gather*}
  |v|_{k,p,\Omega} := 
  \biggl(\sum_{|\delta|=k} |\partial^\delta v|_{0,p,\Omega}^p
   \biggr)^{1/p}, \quad   \|v\|_{k,p,\Omega} := 
  \biggl(\sum_{0 \le m \le k} |v|_{m,p,\Omega}^p \biggr)^{1/p},
\end{gather*}
and $\displaystyle   |v|_{k,\infty,\Omega} := 
  \max_{|\delta|=k} \left\{\mathrm{ess}
   \sup_{\hspace{-5mm}\bfx \in\Omega}|\partial^\delta v(\bfx)|\right\}$,
 $\displaystyle   \|v\|_{k,\infty,\Omega} := 
  \max_{0 \le m  \le k} \left\{|v|_{m,\infty,\Omega}\right\}$.

\subsection{Preliminaries from matrix analysis}
We introduce some facts from the theory of matrix analysis.  For their
proofs, refer to textbooks on matrix analysis such as
\cite{HornJohnson} and \cite{Yamamoto2}.

Let $n \ge 2$ be an integer and $A$ be an $n \times n$ regular matrix.
Note that $A^\top A$ is symmetric positive-definite and has $n$
positive eigenvalues $0 < \mu_1 \le \cdots \le \mu_n$.  
The square roots of $\mu_i$ are called the \textit{singular values} of $A$.
Let $\mu_m :=\mu_1$ and $\mu_M := \mu_n$ be the minimum and maximum eigenvalues.
Then,
\[
  \mu_m |\bfx|^2 \le  |A\bfx|^2 \le \mu_M |\bfx|^2, \quad
  \mu_M^{-1} |\bfx|^2 \le  |A^{-1}\bfx|^2  \le \mu_m^{-1}
   |\bfx|^2,  \quad
    \forall \bfx \in \R^n.
\]
For $A$, the matrix norm $\|A\|$ with respect to the Euclidean norm is 
defined by
\begin{align*}
  \|A\| := \sup_{\bfx \in \R^n}\frac{|A\bfx|}{|\bfx|}.
\end{align*}
From these definitions, we realize that
$\|A\| = \mu_M^{1/2}$ and $\|A^{-1}\| = \mu_m^{-1/2}$.

For the Kronecker product of matrices, we have the following lemma
whose proof is straightforward (see the textbooks mentioned above).
\begin{center}
\fbox{
\begin{minipage}{15truecm}
\begin{lemma}\label{Kronecker}
Let $A$, $B$, $C$, and $D$ be $n\times n$ matrices.  
Then, the following equations hold:
\begin{align*} 
   (A \otimes B)(C \otimes D)  = (AC \otimes BD), \qquad
   (A \otimes B)^\top = A^\top \otimes B^\top.
\end{align*}
Furthermore, if $A$ and $B$ have eigenvalues 
$\lam_i$ and $\mu_j$, $i, j = 1, \cdots, n$, respectively, 
then $\lam_i\mu_j$ are eigenvalues of $A\otimes B$.
\end{lemma}
\end{minipage}
}
\end{center}

\vspace{0.3truecm}
\noindent
\textbf{Exercise:} Prove Lemma~\ref{Kronecker}.

\vspace{0.3truecm}
From Lemma~\ref{Kronecker}, we realize that the minimum and
maximum eigenvalues of
$(A^\top A)\otimes(A^\top A) = (A\otimes A)^\top(A\otimes A)$ are
$0 < \mu_m^2 \le \mu_M^2$.  Hence, for any $\bfw\in \R^{n^2}$,
\begin{align*}
   \mu_m^2 |\bfw|^2 \le |(A \otimes A)\bfw|^2
   & \le \mu_M^2 |\bfw|^2,  \quad
  \mu_M^{-2} |\bfw|^2 \le |(A^{-1}\otimes A^{-1})\bfw|^2
   \le \mu_m^{-2} |\bfw|^2.
\end{align*}

The above facts can be extended straightforwardly to the case of
the higher-order Kronecker product $A\otimes ... \otimes A$.
For $A\otimes ... \otimes A$, $A^{-1}\otimes ... \otimes A^{-1}$
(the $k$th Kronecker products), and we have, for $\bfw \in \R^{n^k}$, 
\begin{align*}
   \mu_m^k |\bfw|^2 & \le |(A \otimes ... \otimes  A)\bfw|^2
   \le \mu_M^k |\bfw|^2, \\
  \mu_M^{-k} |\bfw|^2 & \le |(A^{-1}\otimes ... \otimes A^{-1})\bfw|^2
   \le \mu_m^{-k} |\bfw|^2.
\end{align*}
These inequalities imply that
\begin{align*}
  \|A \otimes ... \otimes  A\| = \|A\|^k, \qquad
 \|A^{-1} \otimes ... \otimes  A^{-1}\| = \|A^{-1}\|^k.
\end{align*}

\subsection{Useful inequalities}
For $N$ positive real numbers $U_1, ..., U_N$, the following
inequalities hold:
\begin{gather}
   \sum_{k=1}^N U_k^p \le N^{\tau(p)}
   \left(\sum_{k=1}^N U_k^2\right)^{p/2}, \quad
   \tau(p) := \begin{cases}
         1-p/2, & 1 \le p \le 2 \\
         0,     & 2 \le p < \infty
     \end{cases},  \label{tau}\\
    \left(\sum_{k=1}^N U_k^2\right)^{p/2} \le N^{\gamma(p)}
   \sum_{k=1}^N U_k^p, \quad
   \gamma(p) := \begin{cases}
         0, &  1 \le p \le 2 \\
         p/2 - 1,     & 2 \le p < \infty
     \end{cases}.
    \label{gamma}
\end{gather}

\vspace{0.3truecm}
\noindent
\textbf{Exercise:} Prove the inequalities \eqref{tau} and \eqref{gamma}.

\subsection{The affine transformation defined by a regular matrix}
Let $A$ be an $n \times n$ matrix with det$A > 0$.
We consider the affine transformation $\varphi(\bfx)$ defined by
$\bfy = \varphi(\bfx) := A \bfx + \mathbf{b}$ for
$\bfx = (x_1, \cdots, x_n)^\top$, $\bfy = (y_1, \cdots, y_n)^\top$
with  $\mathbf{b} \in \R^n$.  Suppose that a reference region
$\widehat\Omega\subset \R^n$ is transformed to a domain $\Omega$ by
$\varphi$; $\Omega := \varphi(\widehat\Omega)$.   Then, a function
$v(\bfy)$ defined on $\Omega$ is pulled-back to the function $\hv(\bfx)$
on $\widehat\Omega$ as    $\hv(\bfx) := v(\varphi(\bfx)) = v(\bfy)$.
Then, we have $\nabla_\bfx \hv = (\nabla_\bfy v) A$,
$\nabla_\bfy v = (\nabla_\bfx \hv) A^{-1}$, and
$|\nabla_\bfy v|^2 = |(\nabla_\bfx \hv) A^{-1}|^2
    = (\nabla_\bfx \hv) A^{-1}A^{-\top} (\nabla_\bfx \hv)^\top$.

The Kronecker product $\nabla\otimes\nabla$ of the gradient $\nabla$ is
defined by
\[
  \nabla\otimes\nabla := \left(
   \frac{\partial\ }{\partial x_1}\nabla, ...,
   \frac{\partial\ }{\partial x_n}\nabla\right)
 = \left(
   \frac{\partial^2\ }{\partial x_1^2},
   \frac{\partial^2\ }{\partial x_1\partial x_2}, ...,
   \frac{\partial^2\ }{\partial x_{n-1}\partial x_n},
   \frac{\partial^2\ }{\partial x_n^2}\right).
\]
We regard $\nabla\otimes\nabla$ to be a row vector.
From this definition, it follows that
\[
   \sum_{|\delta| = 2} (\partial^\delta v)^2 = \sum_{i,j=1}^n 
  \left(\frac{\partial^2 v}{\partial x_i\partial x_j}\right)^2
  = |(\nabla\otimes\nabla) v|^2
\]
and  $(\nabla_\bfx\otimes\nabla_\bfx) \hv =
   \left((\nabla_\bfy\otimes\nabla_\bfy) v\right) (A\otimes A)$, 
   $(\nabla_\bfy\otimes\nabla_\bfy) v =
   \left((\nabla_\bfx\otimes\nabla_\bfx) \hv \right)
   (A^{-1}\otimes A^{-1})$.
Thus, we have
$\|A\|^{-2} |\nabla_\bfx \hv|^2  \le
   |\nabla_\bfy v|^2 \le \|A^{-1}\|^2 |\nabla_\bfx \hv|^2$ and
\begin{align*}
  \sum_{|\delta| = 2} (\partial_\bfy v)^2 & = 
  |(\nabla_\bfy\otimes\nabla_\bfy) v|^2 \\
 & = \left((\nabla_\bfx\otimes\nabla_\bfx) \hv\right) (A^{-1}\otimes A^{-1})
       (A^{-1}\otimes A^{-1})^\top
    \left((\nabla_\bfx\otimes\nabla_\bfx) \hv\right)^\top  \\
 & = \left((\nabla_\bfx\otimes\nabla_\bfx) \hv\right)
      (A^{-1}A^{-\top}\otimes A^{-1}A^{-\top})
       \left((\nabla_\bfx\otimes\nabla_\bfx) \hv\right)^\top, \\
  \|A\|^{-2} \sum_{|\delta| = 2} (\partial_\bfx^\delta \hv)^2 
 & \le \sum_{|\delta| = 2} (\partial_\bfy^\delta v)^2 
 \le \|A^{-1}\|^2 \sum_{|\delta| = 2} (\partial_\bfx^\delta \hv)^2.
\end{align*}
The above inequalities can be easily extended to higher-order
derivatives, and we obtain the following inequalities:
for $k \ge 1$,
\begin{align}
 & \|A\|^{-2k}\sum_{|\delta| = k} (\partial_\bfx^\delta \hv)^2 
 \le \sum_{|\delta| = k} (\partial_\bfy^\delta v)^2 
 \le \|A^{-1}\|^{2k} \sum_{|\delta| = k}
   (\partial_\bfx^\delta \hv)^2,   \notag \\
 & |\!\det{A}|^{1/2} \|A\|^{-k} |\hv|_{k,2,\widehat{\Omega}}
   \le |v|_{k,2,\Omega} \le
   |\!\det{A}|^{1/2} \|A^{-1}\|^k |\hv|_{k,2,\widehat{\Omega}}.
  \label{generalm}
\end{align}

Using the inequalities \eqref{tau} and \eqref{gamma}, we can extend
\eqref{generalm} for the case of arbitrary $p$, $1 \le p < \infty$:
{\allowdisplaybreaks
\begin{align*}
  |v|_{k,p,\Omega}^p & = \int_\Omega \sum_{|\delta|=k}
            |\partial_\bfy^\delta v(\bfy)|^p \dd \bfy
    \le n^{k \tau(p)} \int_\Omega \left(
      \sum_{|\delta|=k} |\partial_\bfy^\delta v(\bfy)|^2 \right)^{p/2}
      \dd \bfy \notag \\
   & \le n^{k\tau(p)} \|A^{-1}\|^{kp}
      \int_\Omega \left(
      \sum_{|\delta|=k} |\partial_\bfx^\delta \hv(\bfx)|^2 \right)^{p/2}
     \dd \bfy \notag \\
   & = n^{k\tau(p)} |\!\det A|\|A^{-1}\|^{kp}
      \int_{\widehat{\Omega}} \left(
      \sum_{|\delta|=k} |\partial_\bfx^\delta \hv(\bfx)|^2 \right)^{p/2}
    \dd \bfx \notag \\
   & \le n^{k(\tau(p)+\gamma(p))} |\!\det A|\|A^{-1}\|^{kp}
     \int_{\widehat{\Omega}} 
      \sum_{|\delta|=k} |\partial_\bfx^\delta \hv(\bfx)|^p  \dd \bfx \notag\\
    & = n^{k (\tau(p)+\gamma(p))} |\!\det A|\|A^{-1}\|^{kp}
   |\hv|_{k,p,\widehat{\Omega}}^p
\end{align*}
}
and
{\allowdisplaybreaks
\begin{align*}
  |v|_{k,p,\Omega}^p & = \int_\Omega \sum_{|\delta|=k}
            |\partial_\bfy^\delta v(\bfy)|^p \dd \bfy \notag 
%
  \ge n^{-k \gamma(p)} \int_\Omega \left(
      \sum_{|\delta|=k} |\partial_\bfy^\delta v(\bfy)|^2 \right)^{p/2}
  \dd \bfy \notag\\
   & \ge n^{-k\gamma(p)} |\!\det A|\|A\|^{-kp}
      \int_\Omega \left(
      \sum_{|\delta|=k} |\partial_\bfx^\delta \hv(\bfx)|^2 \right)^{p/2}
  \dd \bfy \notag\\
   & = n^{ -k\gamma(p)} |\!\det A|\|A\|^{-kp}
      \int_{\widehat{\Omega}} \left(
      \sum_{|\delta|=k} |\partial_\bfx^\delta \hv(\bfx)|^2 \right)^{p/2}
   \dd \bfx \notag\\
   & \ge n^{-k (\tau(p)+\gamma(p))} |\!\det A|\|A\|^{-kp}
      \int_{\widehat{\Omega}} 
      \sum_{|\delta|=k} |\partial_\bfx^\delta \hv(\bfx)|^p  \dd \bfx\notag \\
    & = n^{-k(\tau(p)+\gamma(p))} |\!\det A|\|A\|^{-kp}
    |\hv|_{k,p,\widehat{\Omega}}^p,
\end{align*}
}
where we use the fact that $|v|_{k,2,\Omega}$ contains $n^k$ terms.
Therefore, we obtain the following lemma:
\begin{center}
\fbox{
\begin{minipage}{15truecm}
\begin{lemma}
In the above setting of the linear transformation, we have
\begin{align}
 n^{-k\mu(p)} |\!\det{A}|^{1/p} \|A\|^{-k} |\hv|_{k,p,\widehat{\Omega}}
   \le |v|_{k,p,\Omega} \le n^{k\mu(p)}
   |\!\det{A}|^{1/p} \|A^{-1}\|^k |\hv|_{k,p,\widehat{\Omega}}. 
   \label{generalmp}
\end{align}
where
\[
   \mu(p) := \frac{\tau(p) + \gamma(p)}{p} = 
   \begin{cases}
         1/p -1/2, & 1 \le p \le 2 \\
         1/2 - 1/p,     & 2 \le p \le \infty
     \end{cases}.
\]
\end{lemma}
\end{minipage}
}
\end{center}
\noindent
\textbf{Proof:} We only need to prove the case of $p = \infty$,
and it is done just by letting $p \to \infty$ in 
\eqref{generalmp}.  $\square$

\vvskip
Let us apply \eqref{generalmp} to the case $A \in O(n)$, where
$O(n)$ is the set of orthogonal matrices.  That is,
$A^\top A = AA^\top = I_n$.   In this case, 
$|\mathrm{det}A| = \|A\| = \|A^{-1}\| = 1$.
Thus, we have
\begin{align}
 n^{-k\mu(p)} |\hv|_{k,p,\widehat{\Omega}}
   \le |v|_{k,2,\Omega} \le n^{k\mu(p)}
   |\hv|_{k,p,\widehat{\Omega}}. 
  \label{rotationp}
\end{align}
Those inequalities mean that, if $p = 2$, the Sobolev norms
$|v|_{k,2,\Omega}$ are not affected by rotations.  If $p \neq 2$, 
however, they are affected by rotations up to the constants
$n^{-k\mu(p)}$ and $n^{k\mu(p)}$.

\subsection{The Sobolev imbedding theorem}\label{Sobolevimbedding}
If $1 < p < \infty$, Sobolev's imbedding theorem and Morrey's inequality
imply that
\begin{gather*}
   W^{2,p}(K) \subset C^{1,1-2/p}(K), \quad  p > 2, \\
   H^{2}(K) \subset W^{1,q}(K) \subset C^{0,1-2/q}(K),
     \quad \forall q > 2, \\ 
   W^{2,p}(K) \subset W^{1,2p/(2-p)}(K) \subset C^{0,2(p-1)/p}(K),
    \quad  1 < p < 2.
\end{gather*}
For proofs of the Sobolev imbedding theorems, see
\cite{AdamsFournier} and \cite{Brezis}.
For the case $p=1$, we still have the continuous imbedding
   $W^{2,1}(K) \subset C^0(K)$.
For proof of the critical imbedding, see
\cite[Theorem~4.12]{AdamsFournier} and \cite[Lemma~4.3.4]{BrennerScott}.

\subsection{Gagliardo--Nirenberg's inequality}
\begin{center}
\fbox{
\begin{minipage}{15truecm}
\begin{theorem}[Gagliardo--Nirenberg's inequality]
\label{GagliardoNirenberg}
Let $1 \le p \le \infty$. Let $k$, $m$ be integers such that $k \ge 2$
Then, for $\alpha := m/k$, $0 < \alpha < 1$, 
the following inequality holds:
\begin{equation*}
   |v|_{m,p,\R^n} \le C |v|_{0,p,\R^n}^{1-\alpha} |v|_{k,p,\R^n}^\alpha,
   \quad \forall v \in W^{k,p}(\R^n),
\end{equation*}
where the constant $C$ depends only on $k$, $m$, $p$, and $n$.
\end{theorem}
\end{minipage}
}
\end{center}

For the proof and the general cases of Galliardo--Nirenberg's
inequality, see \cite{Brezis} and the references therein.

\subsection{A standard error analysis of Lagrange interpolation}
In this subsection, we explain a standard error analysis
of Lagrange interpolation.  
First, we prepare a theorem from Ciarlet\cite{Ciarlet}.
Let $\Omega \subset \R^n$ be a bounded domain with the Lipschitz
boundary $\partial\Omega$.  Let $k$ be a positive integer
and $p$ be a real with $1 \le p \le \infty$.
We consider the quotient space $W^{k+1,p}(\Omega)/\PP_k(\Omega)$. As
usual, we introduce the following norm to the space:
\begin{align}
  \|\dot{v}\|_{k+1,p,\Omega} & := \inf_{q \in \PP_k(\Omega)}
  \|v + q\|_{k+1,p,\Omega}, \quad
   \forall \dot{v} \in W^{k+1,p}(\Omega)/\PP_k(\Omega), \\
  \dot{v} & := \left\{ w \in W^{k+1,p}(\Omega)\; |
       \; w - v \in \PP_k(\Omega) \right\}.  \notag
\end{align}
We also define the seminorm of the space by
$|\dot{v}|_{k+1,p,\Omega}:=|v|_{k+1,p,\Omega}$.
Take an arbitrary $q \in \PP_k(\Omega)$.
If $1 \le p < \infty$,  we have
\begin{align*}
   \|v + q\|_{k+1,p,\Omega}^p = |v|_{k+1,p,\Omega}^p
   + \|v + q\|_{k,p,\Omega}^p \ge |v|_{k+1,p,\Omega}^p,
\end{align*}
and if $p = \infty$, we have
\begin{align*}
   \|v + q\|_{k+1,\infty,\Omega} = \max\left\{|v|_{k+1,\infty,\Omega},
   \|v + q\|_{k,\infty,\Omega}\right\} \ge |v|_{k+1,\infty,\Omega}.
\end{align*}
Thus the following inequality follows:
\begin{align}
  |\dot{v}|_{k+1,p,\Omega} \le \|\dot{v}\|_{k+1,p,\Omega}, \qquad
  \forall \dot{v} \in W^{k+1,p}(\Omega)/\PP_k(\Omega).
\end{align}

The next theorem claims the seminorm is actually a norm
of $W^{k+1,p}(\Omega)/\PP_k(\Omega)$.
\begin{center}
\fbox{
\begin{minipage}{15truecm}
\begin{theorem}[Ciarlet\cite{Ciarlet}, Theorem~3.1.1]
There exists a positive constant $C(\Omega)$ depending only on
$k$, $p \in [1,\infty]$, and $\Omega$,  such that the following
estimations hold$:$
\begin{align}
   \|\dot{v}\|_{k+1,p,\Omega} & \le 
  C(\Omega) |\dot{v}|_{k+1,p,\Omega}, \qquad \forall\dot{v}
   \in W^{k+1,p}(\Omega)/\PP_k(\Omega), \notag \\
    \inf_{q \in \PP_k(\Omega)} \|v+q\|_{k+1,p,\Omega} & \le 
  C(\Omega) |v|_{k+1,p,\Omega}, \qquad \forall v
   \in W^{k+1,p}(\Omega). \label{eq3.1.9}
\end{align}
\end{theorem}
\end{minipage}
}
\end{center}
\textbf{Proof:} Let $N$ be the dimension of $\PP_k(\Omega)$ as a vector
space, and $\{q_i\}_{i=1}^N$ be its basis and $\{f_i\}_{i=1}^N$ be the
dual basis of $\{q_i\}$.  That is,
$f_i \in \mathcal{L}(\PP_k(\Omega),\R)$ and they satisfy
$f_i(q_j)=\delta_{ij}$, $i,j=1,\cdots,N$ ($\delta_{ij}$ are Kronecker's
deltas). By Hahn-Banach's theorem, $f_i$ is extended to
$f_i \in \mathcal{L}(W^{k+1,p}(\Omega),\R)$.
For $q \in \PP_k(\Omega)$, we have
\begin{align*}
   q = 0  \Longleftrightarrow f_i(q) = 0, \quad 1 \le i \le N.
\end{align*}

Now, we claim that there exists a constant $C(\Omega)$ such that
\begin{align}
  \|v\|_{k+1,p,\Omega} \le C(\Omega)
  \left(|v|_{k+1,p,\Omega} + \sum_{i=1}^N |f_i(v)|\right), \quad
   \forall v \in W^{k+1,p}(\Omega).
   \label{eq3.1.11}
\end{align}
Suppose that \eqref{eq3.1.11} holds.  For given
$v \in W^{k+1,p}(\Omega)$, let $q \in \PP_k(\Omega)$ be defined with the
extended $f_i \in \mathcal{L}(W^{k+1,p}(\Omega),\R)$ by
\begin{align*}
  q = \sum_{i=1}^N\lam_i q_i, \quad
    \lambda_i := - f_i(v), \quad i = 1, \cdots ,N.
\end{align*}
Then, we have $f_i(v+q)=0$, $i=1,\cdots,N$.  Therefore,
The inequality \eqref{eq3.1.9} follows from \eqref{eq3.1.11}.

We now show the inequality \eqref{eq3.1.11} by contradiction.
Assume that \eqref{eq3.1.11} does not hold.  Then, there exists a
sequence $\{v_l\}_{l=1}^\infty \subset W^{k+1,p}(\Omega)$ such that
\begin{equation*}
   \|v_l\|_{k+1,p,\Omega} = 1, \; \forall l,\qquad
  \lim_{l \to \infty} \left(
   |v_l|_{k+1,p,\Omega} + \sum_{i=1}^N |f_i(v_l)|\right) = 0.
\end{equation*}
By the compactness of the inclusion
$W^{k+1,p}(\Omega) \subset W^{k,p}(\Omega)$, there exists a subsequence
$\{v_{l_m}\}$ and $v \in W^{k,p}(\Omega)$ such that
\[
   \lim_{l_m \to \infty} \|v_{l_m} - v\|_{k,p,\Omega} = 0, \qquad
   \lim_{l_m \to \infty} |v_{l_m}|_{k+1,p,\Omega} = 0.
\]
Here, $\{v_{l_m}\}$ is a Cauchy sequence in $W^{k,p}(\Omega)$.
We show that it is also a Cauchy sequence in $W^{k+1,p}(\Omega)$ as well.
If, for example, $1 \le p < \infty$, we have
\begin{align*}
   \lim_{l_m, l_n \to \infty}  \|v_{l_m} - v_{l_n}\|_{k+1,p,\Omega} 
    & = \lim_{l_m, l_n \to \infty} \left(
    \|v_{l_m} - v_{l_n}\|_{k,p,\Omega}^p
   + |v_{l_m} - v_{l_n}|_{k+1,p,\Omega}^p\right)^{1/p} \\
  & \le  \lim_{l_m, l_n \to \infty} \left(
    \|v_{l_m} - v_{l_n}\|_{k,p,\Omega}^p
   + 2^{p-1}(|v_{l_m}|_{k+1,p,\Omega}^p + |v_{l_n}
  |_{k+1,p,\Omega}^p)\right)^{1/p}\\
  &  = 0.
\end{align*}
The case for $p = \infty$ is similarly shown.
Hence, $v$ belong s to $W^{k+1,p}(\Omega)$, and $\{v_{l_m}\}$ satisfies
\[
   \lim_{l_m \to \infty} \|v_{l_m} - v\|_{k+1,p,\Omega} = 0.
\]
This $v \in W^{k+1,p}(\Omega)$ satisfies 
\begin{align*}
   |\partial^\beta v|_{0,p,\Omega} = \lim_{l_k \to 0}
     |\partial^\beta v_{l_k}|_{0,p,\Omega} = 0, \quad
    \forall \beta, \quad |\beta| = k+1,
\end{align*}
and thus $v \in \PP_k(\Omega)$.  Therefore, because 
\[
   \sum_{i=1}^N |f_i(v)| = \lim_{l_m \to \infty}
   \sum_{i=1}^N |f_i(v_{l_m})| =0,
\]
we conclude $v = 0$.   However, this contradicts to
$\displaystyle
\|v\|_{k+1,p,\Omega}=\lim_{l_m \to \infty}\|v_{l_m}\|_{k+1,p,\Omega}=1$.
 $\square$

\vspace{2mm}
We are now ready to prove the first inequality in Theorem~\ref{Thm314}.
Recall that $\hK$ is the reference triangle and $K$ is mapped
as $K = \varphi(\hK)$ with $\varphi(\bfx) = A \bfx + \bfb$.
\begin{center}
\fbox{
\begin{minipage}{15truecm}
\begin{theorem}
Suppose that $\|A^{-1}\| \ge 1$.  Then, there exists a constant \\
$C =  C(\hK,p,k,m)$ independent of $K$ such that
\begin{align}
  \|v - \I_K^k v\|_{m,p,K} \le C \|A\|^{k+1}\|A^{-1}\|^m
   |v|_{k+1,p,K}, \quad \forall v \in W^{k+1,p}(K).
   \label{standard-est}
\end{align}
\end{theorem}
\end{minipage}
}
\end{center}
\noindent
\textbf{Proof:} Note that, for arbitrary
$\hat{v} \in W^{k+1,p}(\hK)$ and $\hat{p} \in \PP_k(\hK))$, we have
\begin{align*}
   \hat{v} - \I_{\hK}^k \hat{v} = (I - \I_{\hK}^k \hat{v})
   (\hat{v} + \hat{p}),
\end{align*}
where $I: W^{k+1,p}(\hK) \to W^{m,p}(\hK)$ is the identity mapping,
which is obviously continuous.  Therefore, it follows from
\eqref{eq3.1.9} that
\begin{align*}
   \left\|\hat{v} - \I_{\hK}^k \hat{v}\right\|_{m,p,\hK} & \le
   \left\|I - \I_{\hK}^k\right\|_{\LL(W^{k+1,p}(\hK),W^{m,p}(\hK))}
   \inf_{\hat{p} \in \PP_k(\hK)}
   \left\|\hat{v} + \hat{p}\right\|_{k+1,p,\hK} \\
   & \le C_1|\hat{v}|_{k+1,p,\hK},
\end{align*}
where the constant $C_1$ depends on $\hK$, $m$, $k$, $p$,
(and $\I_{\hK}^k$).

Note that the mapping between $W^{n,p}(K)$ and $W^{n,p}(\hK)$
($n = m$ or $n = k+1$) defined by the pull-back
$\hat{v} = v\circ \varphi$ is an isomorphism.  By \eqref{generalmp},
we have
\begin{align*}
   \left\|v - \I_{K}^k v\right\|_{m,p,K} & = \left(
   \sum_{l = 0}^m \left|v - \I_{K}^k v\right|_{l,p,K}^p
  \right)^{1/p} \\
 & \le \left(\sum_{l=0}^m
  n^{l\mu(p)}|\det A|\|A^{-1}\|^{lp}
   \left|\hat{v} - \I_{\hK}^k \hat{v}\right|_{l,p,\hK}^p
  \right)^{1/p}, \\
 & \le n^{m\mu(p)}|\det A|^{1/p}\|A^{-1}\|^m 
   \left\|\hat{v} - \I_{\hK}^k \hat{v}\right\|_{m,p,\hK}, \\
  \left|\hat{v}\right|_{k+1,p,K} & \le
  n^{(k+1)\mu(p)}|\det A|^{-1/p}\|A\|^{k+1}
   \left|v\right|_{k+1,p,\hK},
\end{align*}
because of the assumption $\|A^{-1}\| \ge 1$.
Combining these inequalities, the proof is completed with
$C := n^{(k+1+m)\mu(p)} C_1$.  $\square$

\vspace{3mm}
Combining these propositions with Lemma~\ref{inscribed}, we see that,
for arbitrary $v \in W^{k+1,p}(K)$,
\begin{align*}
   \|v - \I_K^k v\|_{m,p,K} & \le C \|A\|^{k+1}\|A^{-1}\|^m
   |v|_{k+1,p,K} 
   \le C \left(\frac{h_K}{\rho_{\hK}}\right)^{k+1}
    \left(\frac{h_{\hK}}{\rho_K}\right)^m |v|_{k+1,p,K} \\
  & \le C \frac{h_{\hK}^m}{\rho_{\hK}^{k+1}}
    \frac{h_K^{k+1}}{\rho_K^m} |v|_{k+1,p,K}.
\end{align*}
If there exists a constant $\sigma$ such that
$h_K/\rho_K \le \sigma$, then $\rho_K^{-1} \le \sigma h_K^{-1}$, and
we obtain the following standard error estimation.

\begin{center}
\fbox{
\begin{minipage}{15truecm}
\begin{theorem}
Let $K \subset \R^2$ be a triangle with $h_K \le 1$.
Suppose that $h_K/\rho_K \le \sigma$, where $\sigma$ is a positive
constant.  Then, there exists a constant $C =  C(\hK,p,k,m, \sigma)$
independent of $K$ such that
\begin{align}
  \|v - \I_K^k v\|_{m,p,K} \le C h_K^{k+1-m}
   |v|_{k+1,p,K}, \quad \forall v \in W^{k+1,p}(K).
   \label{standard-est2}
\end{align}
\end{theorem}
\end{minipage}
}
\end{center}

\section{Babu\v{s}ka--Aziz's technique}
In the previous section, we have proved the standard
error estimation \eqref{standard-est}, \eqref{standard-est2}. 
To improve them, we introduce the technique given by
Babu\v{s}ka--Aziz \cite{BabuskaAziz}. 

Let $\hK$ be the reference triangle with the vertices $(0,0)^\top$,
$(1,0)^\top$, and $(0,1)^\top$.   For $\hK$, the sets
$\Xi_p^i \subset W^{1,p}(\hK)$, $i = 1, 2$, $p \in [1,\infty]$ are
defined by
\begin{align*}
  \Xi_p^{(1,0),1} & := \left\{ v \in W^{1,p}(\hK) \Bigm| 
    \int_0^1 v(s,0) \dd s = 0 \right\}, \\
  \Xi_p^{(0,1),1} & := \left\{ v \in W^{1,p}(\hK) \Bigm| 
    \int_0^1 v(0,s) \dd s = 0 \right\}.
\end{align*}
The constant $A_p$ is then defined by
\begin{align*}
   A_p := \sup_{v \in \Xi_p^{(1,0),1}} \frac{|v|_{0,p,\hK}}
  {|v|_{1,p,\hK}} =
  \sup_{v \in \Xi_p^{(0,1),1}} \frac{|v|_{0,p,\hK}}
  {|v|_{1,p,\hK}}, \qquad 1 \le p \le \infty.
\end{align*}
The second equation in the above definition follows from the symmetry of
$\hK$.  The constant $A_p$ (and its reciprocal $1/A_p$) is called the
\textbf{Babu\v{s}ka--Aziz constant} for $p \in [1,\infty]$.  According to
Liu--Kikuchi \cite{LiuKikuchi}, $A_2$ is the maximum positive solution of
the equation $1/x + \tan(1/x) = 0$, and $A_2 \approx 0.49291$.

In the following, we show that $A_p < \infty$
(Babu\v{s}ka--Aziz \cite[Lemma~2.1]{BabuskaAziz} and
Kobayashi--Tsuchiya  \cite[Lemma~1]{KobayashiTsuchiya1}).

\begin{center}
\fbox{
\begin{minipage}{15truecm}
\begin{lemma}\label{L1.2} We have $A_p < \infty$, $p \in [1,\infty]$.
\end{lemma}
\end{minipage}
}
\end{center}
\textbf{Proof:} The proof is by contradiction.  Assume that $A_p = \infty$.
Then, there exists a sequence $\{u_k\}_{i=1}^\infty \subset \Xi_p^{(1,0),1}$
such that
\[
   |u_k|_{0,p,\hK} = 1, \qquad \lim_{k \to \infty} |u_k|_{1,p,\hK} = 0.
\]
From the inequality \eqref{eq3.1.9}, for an arbitrary $\varepsilon > 0$,
there exists a sequence $\{q_k\}\subset \PP_0(\hK)$ such that
\begin{gather*}
   \inf_{q \in \PP_0(\hK)}
  \|u_k + q\|_{1,p,\hK} 
   \le \|u_k + q_k\|_{1,p,\hK}
   \le \inf_{q \in \PP_0(\hK)}
  \|u_k + q\|_{1,p,\hK} + \frac{\varepsilon}{k}
   \le C |u_k|_{1,p,\hK} + \frac{\varepsilon}{k}, \\
  \lim_{k \to \infty} \|u_k + q_k\|_{1,p,\hK} = 0.
\end{gather*}
Since the sequence $\{u_k\}\subset W^{1,p}(\hK)$ is bounded,
$\{q_k\} \subset \PP_0(\hK) = \R$ is also bounded. Therefore, there exists a
subsequence $\{q_{k_i}\}$ such that $q_{k_i}$ converges to
$\bar{q} \in \PP_0(\hK)$.  Thus, in particular, we have
\[
   \lim_{k_i \to \infty} \|u_{k_i} + \bar{q}\|_{1,p,\hK} = 0.
\]
Let $\Gamma$ be the edge of $\hK$ connecting $(1,0)^\top$ and
$(0,0)^\top$ and $\gamma : W^{1,p}(\hK) \to W^{1-1/p,p}(\Gamma)$ be the
trace operator.  The continuity of $\gamma$
and the inclusion $W^{1-1/p,p}(\Gamma) \subset L^1(\Gamma)$ yield
\[
 0 = \lim_{k_i \to \infty} \int_\Gamma \gamma(u_{k_i} + \bar{q}) \dd s
  = \int_\Gamma \bar{q} \dd s,
\]
because $u_{k_i} \in \Xi_p^1$.  Thus, we find that $\bar{q} = 0$ and
$\lim_{k_i \to \infty}\|u_{k_i}\|_{1,p,\hK} = 0$.  This contradicts
$\lim_{k_i \to \infty}\|u_{k_i}\|_{1,p,\hK} \ge
 \lim_{k_i \to \infty}|u_{k_i}|_{0,p,\hK}=1$.
$\square$

\vspace{0.4truecm}
We define the bijective linear transformation $\Fab:\R^2 \to \R^2$ by
\[
   (x^*, y^*)^\top = (\alpha x, \beta y)^\top, \qquad
   (x,y)^\top \in \R^2, \quad \alpha, \beta > 0.
\]
The map $\Fab$ is called the \textbf{squeezing transformation}.

Now, we consider the ``squeezed'' triangle  $\Kab:= \Fab(\hK)$.
Take an arbitrary $v \in W^{2,p}(\Kab)$, and pull-back $v$ to
$u := v \circ \Fab \in W^{2,p}(\hK)$.  For, $p$, $1 \le p < \infty$,
we have
{\allowdisplaybreaks
\begin{gather}
  \frac{|v|_{0,p,\Kab}^p}{|v|_{1,p,\Kab}^p}
    = \frac{|u|_{0,p,\hK}^p}{\frac{1}{\alpha^p}|u_{x}|_{0,p,\hK}^p
                + \frac{1}{\beta^p}|u_{y}|_{0,p,\hK}^p},  \label{trans1} \\
  \frac{|v|_{0,p,\Kab}^p}{|v|_{2,p,\Kab}^p}
    = \frac{|u|_{0,p,\hK}^p}
     {\frac{1}{\alpha^{2p}}|u_{xx}|_{0,p,\hK}^p 
    + \frac{2}{\alpha^p\beta^p} |u_{xy}|_{0,p,\hK}^p 
        + \frac{1}{\beta^{2p}} |u_{yy}|_{0,p,\hK}^p},  \label{trans2} \\
  \frac{|v|_{1,p,\Kab}^p}{|v|_{2,p,\Kab}^p}
    = \frac{\frac{1}{\alpha^p}|u_x|_{0,p,\hK}^p
            + \frac{1}{\beta^p}|u_y|_{0,p,\hK}^p}
       {\frac{1}{\alpha^{2p}}|u_{xx}|_{0,p,\hK}^p 
          + \frac{2}{\alpha^p\beta^p}|u_{xy}|_{0,p,\hK}^p
       + \frac{1}{\beta^{2p}} |u_{yy}|_{0,p,\hK}^p}.        
     \label{trans3}
\end{gather}
}
In the following we explain how these equations are derived.

Note that, for
$(x,y)^\top \in \hK$ and $(x^*,y^*)^\top = (\alpha x,\beta y)\top\in \Kab$,
we have
\[
   u_{x} = \alpha v_{x^*}, \qquad u_{y} = \beta v_{y^*}, \quad
\]
and
\begin{align*}
 |v_{x^*}|_{0,p,\Kab}^p = \int_{\Kab} |v_{x^*}|^p \dd\bfx^*
   =  \frac{1}{\alpha^p}\int_{\Kab} |u_{x}|^p \dd \bfx^*
   = \frac{\beta}{\alpha^{p-1}} \int_{\hK} |u_{x}|^p \dd \bfx
   = \frac{\beta}{\alpha^{p-1}} |u_{x}|_{0,p,\hK}^p.
\end{align*}
Here, $\dd {\bfx} := \dd x \dd y$, 
$\dd {\bfx^*} := \dd x^* \dd y^*$, and used the fact
$\det(D\Fab) = \alpha\beta$,
where $D \Fab$ is the Jacobian matrix of $\Fab$.
Similarly, we obtain
\begin{align*}
 |v_{y^*}|_{0,p,\Kab}^p = \frac{\alpha}{\beta^{p-1}}
   |u_{y}|_{0,p,\hK}^p, \qquad
  |v|_{0,p,\Kab}^p = \alpha\beta |u|_{0,p,\hK}^p.
\end{align*}
Therefore, these equations yield \eqref{trans1}:
{\allowdisplaybreaks
\begin{align*} 
  \frac{|v|_{0,p,\Kab}^p}{|v|_{1,p,\Kab}^p}
    = \frac{\alpha\beta |u|_{0,p,\hK}^p}
      {\frac{\beta}{\alpha^{p-1}} |u_{x}|_{0,p,\hK}^p
        + \frac{\alpha}{\beta^{p-1}} |u_{y}|_{0,p,\hK}^p}
    =  \frac{|u|_{0,p,\hK}^p} 
      {\frac{1}{\alpha^{p}}|u_{x}|_{0,p,\hK}^p
     + \frac{1}{\beta^{p}}|u_{y}|_{0,p,\hK}^p}.
\end{align*}
}
Similarly, the equations
\begin{align*}
   |v_{x^*x^*}|_{0,p,\Kab}^p & = \frac{\beta}{\alpha^{2p-1}}|u_{xx}|_{0,p,\hK}^p, \\
   |v_{x^*y^*}|_{0,p,\Kab}^p & = \frac{1}{\alpha^{p-1}\beta^{p-1}}
   |u_{xy}|_{0,p,\hK}^p, \\
  |v_{y^*y^*}|_{0,p,\Kab}^p & = \frac{\alpha}{\beta^{2p-1}}
   |u_{yy}|_{0,p,\hK}^p
\end{align*}
are obtained and yield \eqref{trans2} and \eqref{trans3} as
{\allowdisplaybreaks
\begin{align*}
 \frac{|v|_{0,p,\Kab}^p}{|v|_{2,p,\Kab}^p}
  & = \frac{\alpha\beta |u|_{0,p,\hK}^p}
       {\frac{\beta}{\alpha^{2p-1}} |u_{xx}|_{0,p,\hK}^p
           + \frac{2}{\alpha^{p-1}\beta^{p-1}} |u_{xy}|_{0,p,\hK}^p
           + \frac{\alpha}{\beta^{2p-1}} |u_{yy}|_{0,p,\hK}^p} \\
  &   = \frac{|u|_{0,p,\hK}^p}
        {\frac{1}{\alpha^{2p}}|u_{xx}|_{0,p,\hK}^p
           + \frac{2}{\alpha^p\beta^p} |u_{xy}|_{0,p,\hK}^p
           + \frac{1}{\beta^{2p}} |u_{yy}|_{0,p,\hK}^p}, \\
 \frac{|v|_{1,p,\Kab}^p}{|v|_{2,p,\Kab}^p}
  &  = \frac{\frac{\beta}{\alpha^{p-1}} |u_x|_{0,p,\hK}^p
           + \frac{\alpha}{\beta^{p-1}} |u_y|_{0,p,\hK}^p}
      {\frac{\beta}{\alpha^{2p-1}} |u_{xx}|_{0,p,\hK}^p
         + \frac{\alpha}{\beta^{2p-1}} |u_{yy}|_{0,p,\hK}^p
         + \frac{2}{\alpha^{p-1}\beta^{p-1}} |u_{xy}|_{0,p,\hK}^p} \\
 &   = \frac{\frac{1}{\alpha^p}|u_x|_{0,p,\hK}^p
              + \frac{1}{\beta^p}|u_y|_{0,p,\hK}^p}
        {\frac{1}{\alpha^{2p}}|u_{xx}|_{0,p,\hK}^p
            + \frac{2}{\alpha^p\beta^p} |u_{xy}|_{0,p,\hK}^p
            + \frac{1}{\beta^{2p}} |u_{yy}|_{0,p,\hK}^p}.
\end{align*}
}

Next, let $p=\infty$.  Then, we have
\begin{align*}
   |v|_{0,\infty,\Kab} & = |u|_{0,\infty,\hK}, \qquad
   |v|_{1,\infty,\Kab} = \max \left\{|u_{x}|_{1,\infty,\hK}/{\alpha},
     |u_{y}|_{1,\infty,\hK}/{\beta}\right\}, \\
   |v|_{2,\infty,\Kab} & = \max \left\{
     |u_{xx}|_{2,\infty,\hK}/{\alpha^2},
    |u_{xy}|_{2,\infty,\hK}/(\alpha\beta),
   |u_{yy}|_{2,\infty,\hK}/{\beta^2} \right\},
\end{align*}
and obtain
\begin{gather}
  \frac{|v|_{0,\infty,\Kab}}{|v|_{1,\infty,\Kab}}
    = \frac{|u|_{0,\infty,\hK}}
       {\max \left\{\frac{1}{\alpha}|u_{x}|_{0,\infty,\hK},
       \frac{1}{\beta}|u_{y}|_{0,\infty,\hK} \right\}},
        \label{trans1-inf} \\
  \frac{|v|_{0,\infty,\Kab}}{|v|_{2,\infty,\Kab}} =
   \frac{|u|_{0,\infty,\hK}}
    {\max \left\{\frac{1}{\alpha^2}|u_{xx}|_{0,\infty,\hK},
          \frac{1}{\alpha\beta}|u_{xy}|_{0,\infty,\hK},
          \frac{1}{\beta^2} |u_{yy}|_{0,\infty,\hK} \right\} },
         \label{trans2-inf} \\
  \frac{|v|_{1,\infty,\Kab}}{|v|_{2,\infty,\Kab}}
    = \frac{\max \left\{\frac{1}{\alpha}|u_x|_{0,\infty,\hK},
         \frac{1}{\beta} |u_y|_{0,\infty,\hK} \right\} }
      {\max \left\{\frac{1}{\alpha^2}|u_{xx}|_{0,\infty,\hK},
       \frac{1}{\alpha\beta} |u_{xy}|_{0,\infty,\hK},
        \frac{1}{\beta^2} |u_{yy}|_{0,\infty,\hK} \right\}}.
    \label{trans3-inf}
\end{gather}

\vspace{0.4truecm}
For a triangle $K$ and $1 \le p \le \infty$, we define
$\T_p^1(K) \subset W^{2,p}(K)$ by
\begin{align*}
   \T_p^1(K) & := \left\{ v \in W^{2,p}(K) \,
  \bigm| \, v(\bfx_i) = 0, \; i = 1,2,3\right\}.
\end{align*}
Note that if $v \in \T_p^1(\Kab)$, then
$u := v\circ \Fab \in \T_p^1(\hK)$.

The following lemma is from
Babu\v{s}ka--Aziz \cite[Lemma~2.2]{BabuskaAziz} and
Kobayashi--Tsuchiya  \cite[Lemma~3]{KobayashiTsuchiya1}.

\begin{center}
\fbox{
\begin{minipage}{15truecm}
\begin{lemma} \label{L2.1}
The constant $B_p^{1,1}(\Kab)$ is defined by
\[
  B_p^{1,1}(\Kab) := \sup_{v \in \T_p^1(\Kab)} 
  \frac{|v|_{1,p,\Kab}}{|v|_{2,p,\Kab}}, \qquad 1 \le p \le \infty.
\]
Then, we have $B_p^{1,1}(\Kab) \le \max\{\alpha,\beta\} A_p$.

\end{lemma}
\end{minipage}
}
\end{center}
\textbf{Proof:} Suppose first that $1 \le p < \infty$.
Take an arbitrary $v \in \T_p^1(\Kab)$ and define
$u \in \T_p^1(\hK)$ by $u(x,y):=v(x^*,y^*)$,
$(x^*,y^*)^\top = (\alpha x,\beta y)^\top$.  By \eqref{trans3}, we find
\begin{align*}
  \frac{|v|_{1,p,\Kab}^p}{|v|_{2,p,\Kab}^p}
   &  = \frac{\frac{1}{\alpha^p} |u_x|_{0,p,\hK}^p
          + \frac{1}{\beta^p}|u_y|_{0,p,\hK}^p }
    {\frac{1}{\alpha^{2p}} |u_{xx}|_{0,p,\hK}^p 
  + \frac{1}{\alpha^p\beta^p}  |u_{xy}|_{0,p,\hK}^p + 
     \frac{1}{\alpha^p\beta^p} |u_{xy}|_{0,p,\hK}^p 
     + \frac{1}{\beta^{2p}} |u_{yy}|_{0,p,\hK}^p } \\
   & \le \frac{\max\{\alpha^p,\beta^p\}
      \left(\frac{1}{\alpha^{p}} |u_x|_{0,p,\hK}^p
       + \frac{1}{\beta^p}|u_y|_{0,p,\hK}^p\right)}
      { \frac{1}{\alpha^{p}} \left(|u_{xx}|_{0,p,\hK}^p
     + |u_{xy}|_{0,p,\hK}^p\right)
        + \frac{1}{\beta^p}\left(
     |u_{xy}|_{0,p,\hK}^p + |u_{yy}|_{0,p,\hK}^p \right)} \\
   & =  \max\{\alpha^p,\beta^p\}
     \frac{\frac{1}{\alpha^p}|u_x|_{0,p,\hK}^p
         + \frac{1}{\beta^p} |u_y|_{0,p,\hK}^p}
      {\frac{1}{\alpha^p}|u_{x}|_{1,p,\hK}^p
        + \frac{1}{\beta^p}|u_{y}|_{1,p,\hK}^p}.
\end{align*}
Here, we used the fact that, for $X$, $Y > 0$,
\begin{align*}
  \frac{1}{\frac{1}{\alpha^p}X + \frac{1}{\beta^p}Y} \le 
  \frac{\max\{\alpha^p,\beta^p\}}{X + Y}.
\end{align*}
Note that $u(0,0) = u(1,0) = 0$ by the definition of $\T_p^1(\hK)$
and $u_x \in \Xi_p^{(1,0),1}$.  Thus, by Lemma~\ref{L1.2}, we realize
that
\begin{align*}
  |u_x|_{0,p,\hK}^p \le A_p^{p} |u_x|_{1,p,\hK}.
\end{align*}
By the same reason, we realize that $u_y \in \Xi_p^{(0,1),1}$ and
\begin{align*}
  |u_y|_{0,p,K}^p \le A_p^{p} |u_y|_{1,p,K}^p.
\end{align*}
Inserting those inequalities into the above estimation, we obtain
\begin{align*}
  \frac{|v|_{1,p,\Kab}^p}{|v|_{2,p,\Kab}^p}
  \le \max\{\alpha^p,\beta^p\} 
    \frac{\frac{A_p^p}{\alpha^p} |u_x|_{1,p,\hK}^p
          + \frac{A_p^p}{\beta^p}|u_y|_{1,p,\hK}^p}
        {\frac{1}{\alpha^p} |u_x|_{1,p,\hK}^p 
       + \frac{1}{\beta^p}|u_y|_{1,p,\hK}^p} =
       \left(\max\{\alpha,\beta\}\right)^p  A_p^p,
\end{align*}
and conclude
\begin{align*}
 B_p^{1,1}(\Kab) = \sup_{v \in \T_p^1(\Kab)}
    \frac{|v|_{1,p,\Kab}}{|v|_{2,p,\Kab}} \le
    \max\{\alpha,\beta\} A_p.
\end{align*}

Next, let $p=\infty$.  By \eqref{trans3-inf}, we immediately obtain
{\allowdisplaybreaks
\begin{align*}
    \frac{|v|_{1,\infty,\Kab}}{|v|_{2,\infty,\Kab}}
   & = \frac{\max \left\{ \frac{1}{\alpha}|u_x|_{0,\infty,\hK},
          \frac{1}{\beta} |u_y|_{0,\infty,\hK} \right\} }
  { \max \left\{ \max \left\{ \frac{|u_{xx}|_{0,\infty,\hK}}{\alpha^{2}}, 
                     \frac{|u_{xy}|_{0,\infty,\hK}}{\alpha\beta}\right\},
     \max \left\{ \frac{|u_{xy}|_{0,\infty,\hK}}{\alpha\beta},
            \frac{|u_{yy}|_{0,\infty,\hK}}{\beta^2}, \right\}\right\}} \\
  & \le  \frac{ \max\{\alpha,\beta\}
    \max \left\{\frac{1}{\alpha} |u_x|_{0,\infty,\hK},
           \frac{1}{\beta} |u_y|_{0,\infty,\hK} \right\} }
   {\max \left\{ \frac{1}{\alpha}\max
   \left\{|u_{xx}|_{0,\infty,\hK}, |u_{xy}|_{0,\infty,\hK}\right\},
    \frac{1}{\beta} 
   \max \left\{ |u_{xy}|_{0,\infty,\hK}, |u_{yy}|_{0,\infty,\hK},
                \right\}  \right\}} \\
  & =  \max\{\alpha,\beta\}
   \frac{\max \left\{ \frac{1}{\alpha}|u_x|_{0,\infty,\hK},
           \frac{1}{\beta} |u_y|_{0,\infty,\hK} \right\} }
   {\max \left\{ \frac{1}{\alpha}|u_{x}|_{1,\infty,\hK},
    \frac{1}{\beta} |u_{y}|_{1,\infty,\hK} \right\}} \\
   & \le \max\{\alpha,\beta\} 
   \frac{A_\infty \max \left\{ \frac{1}{\alpha}|u_x|_{1,\infty,\hK},
               \frac{1}{\beta} |u_y|_{1,\infty,\hK} \right\}}
  {\max \left\{ \frac{1}{\alpha}|u_x|_{1,\infty,\hK},
          \frac{1}{\beta}|u_y|_{1,\infty,\hK} \right\} }
    = \max\{\alpha,\beta\}   A_\infty.
   \hspace{1.6cm}  \square
\end{align*}
}

\vvskip
The following lemma is from
Babu\v{s}ka--Aziz \cite[Lemma~2.3,2.4]{BabuskaAziz} and
Kobayashi--Tsuchiya  \cite[Lemma~4,5]{KobayashiTsuchiya1}.

\begin{center}
\fbox{
\begin{minipage}{15truecm}
\begin{lemma}\label{L2.2} The constants $B_p^{0,1}(\Kab)$,
$\widetilde{A}_p$ are defined by
\[
  B_p^{0,1}(\Kab) := \sup_{v \in \T_p^1(\Kab)} 
  \frac{|v|_{0,p,\Kab}}{|v|_{2,p,\Kab}}, \quad
  \widetilde{A}_p := B_p^{0,1}(\hK) := \sup_{v \in \T_p^1(\hK)} 
  \frac{|v|_{0,p,\hK}}{|v|_{2,p,\hK}}, \; 1 \le p \le \infty.
\]
Then, we have the estimation
$B_p^{0,1}(\Kab) \le \max\{\alpha^2,\beta^2\}
   \widetilde{A}_p < +\infty$.
\end{lemma}
\end{minipage}
}
\end{center}
\textbf{Proof: } The proof of $\widetilde{A}_p < +\infty$ is 
very similar to that of Lemma~\ref{L1.2} and is by contradiction. 
Supposet that $\widetilde{A}_p = \infty$.  Then, there exists 
$\{u_m\}_{m=1}^\infty \subset \T_p^1(\hK)$ such that
\[
   |u_m|_{0,p,\hK} = 1, \qquad
   \lim_{m \to \infty} |u_m|_{2,p,\hK} = 0.
\]
Then, by \eqref{eq3.1.9}, there exists $\{q_m\}\subset \PP_1(\hK)$ such that
\begin{align*}
  \lim_{m \to \infty} & \|u_m + q_m\|_{2,p,\hK} = 0.
\end{align*}
Since $|u_m|_{0,p,\hK}$ and $|u_m|_{2,p,\hK}$ are bounded,
$|u_m|_{1,p,\hK}$ and $\|u_m\|_{2,p,\hK}$ are bounded as well by
Gagliardo--Nirenberg's inequality (Theorem~\ref{GagliardoNirenberg}).  Hence,
 $\{q_m\} \subset \PP_1(\hK)$ is also bounded.  Thus, there exists
a subsequence $\{q_{m_i}\}$ which converges to $\bar{q} \in \PP_1(\hK)$.
In particular, we have
\[
   \lim_{m_i \to \infty} \|u_{m_i} + \bar{q}\|_{2,p,\hK} = 0.
\]
Since $\{u_{m}\} \subset \T_p^1(\hK)$, we conclude that
$\bar{q}\in \T_p^1(\hK) \cap \PP_1(\hK)$ and $\bar{q} = 0$.
Therefore, we reach $\lim_{m_i \to \infty}\|u_{m_i}\|_{2,p,\hK} = 0$
which contradicts to $\displaystyle \lim_{m_i \to \infty}
\|u_{m_i}\|_{2,p,\hK} \ge  \lim_{m_i \to \infty}|u_{m_i}|_{0,p,\hK}=1$.

We now consider the estimation for the case $1 \le p < \infty$.
From \eqref{trans2} we have
\begin{align*}
  \frac{|v|_{0,p,\Kab}^p}{|v|_{2,p,\Kab}^p}
   & = \frac{|u|_{0,p,\hK}^p}
    { \frac{1}{\alpha^{2p}}|u_{xx}|_{0,p,\hK}^p
     + \frac{2}{\alpha^p\beta^p} |u_{xy}|_{0,p,\hK}^p
 + \frac{1}{\beta^{2p}} |u_{yy}|_{0,p,\hK}^p} \\
 & \le \frac{\max\{\alpha^{2p},\beta^{2p}\}|u|_{0,p,\hK}^p}
  {|u_{xx}|_{0,p,\hK}^p + 2|u_{xy}|_{0,p,\hK}^p + |u_{yy}|_{0,p,\hK}^p}
         \le 
   \left(\max\{\alpha^2,\beta^2\}\right)^p\widetilde{A}_p^p,
\end{align*}
and Lemma is shown for this case.
The proof for the case $p=\infty$ is very similar.  $\square$

\vspace{0.3truecm}
\noindent
\textbf{Exercise:} In Lemma~\ref{L2.2}, prove the case $p=\infty$.

\vspace{0.3truecm}
We may apply Lemmas~\ref{L2.1} and \ref{L2.2} to 
$v - \I_{\Kab}^1 v \in \T_p^1(\Kab)$ for
$v \in W^{2,p}(\Kab)$, and obtain the following corollary.
\begin{center}
\fbox{
\begin{minipage}{15truecm}
\begin{corollary}\label{cor2.2}
For arbitrary $v \in W^{2,p}(\Kab)$ $(1 \le p \le \infty)$, the
following estimations hold: 
\begin{align*}
  |v - \I_{\Kab}^1 v|_{1,p,\Kab} & \le 
   \max\{\alpha,\beta\}A_p |v|_{2,p,\Kab}, \\
  |v - \I_{\Kab}^1 v|_{0,p,\Kab} & \le 
   \left(\max\{\alpha,\beta\}\right)^2
   \widetilde{A}_p |v|_{2,p,\Kab}.
\end{align*}
\end{corollary}
\end{minipage}
}
\end{center}

\section{Extending Babu\v{s}ka-Aziz's technique to the higher order
 Lagrange interpolation}
In this section, we prove the following theorem using Babu\v{s}ka-Aziz's
technique.  Let $k$ be a positive integer and $p$ be such that
$1 \le p \le \infty$.  The set $\T_p^k(K)$ is defined by
\begin{align*}
  \T_p^k(K) := \left\{v \in W^{k+1,p}(K) \bigm|
   v(\bfx) = 0, \forall \bfx \in \Sigma^k(K)\right\},
\end{align*}
where $\Sigma^k(K)$ is defined by \eqref{Sigma}.
Note that if $v \in \T_p^k(\Kab)$, then $u = v\circ\Fab \in \T_p^k(\hK)$.

\begin{center}
\fbox{
\begin{minipage}{15truecm}
\begin{theorem}\label{squeezingtheorem}
Take arbitrary $\alpha > 0$ and $\beta > 0$. Then, there exists
a constant $C_{k,m,p}$ such that, for $m = 0, 1, \cdots, k$, 
\begin{equation}
   B_p^{m,k}(\Kab) := \sup_{v \in \T_p^{k}(\Kab)}
  \frac{|v|_{m,p,\Kab}}{|v|_{k+1,p,\Kab}}
   \le \left(\max\{\alpha,\beta\}\right)^{k+1-m}C_{k,m,p}.
  \label{extension}
\end{equation}
Here, $C_{k,m,p}$ depends only on $k$, $m$, and $p$, and
is independent of $\alpha$ and $\beta$.
\end{theorem}
\end{minipage}
}
\end{center}

\begin{figure}[thb]
\begin{center}
  \includegraphics[width=7truecm]{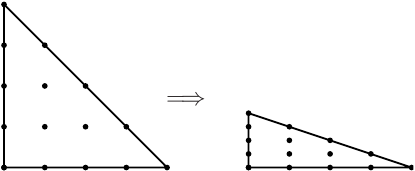}
\caption{Squeezing the reference triangle $\hK$ perpendicularly does not
deteriorate the approximation property of Lagrange interpolation.}
 \label{squeezed}
\end{center}
\end{figure}

\vspace{0.3truecm}
Applying Theorem~\ref{squeezingtheorem} to 
$v - \I_{\Kab}^k v \in \T_p^k(\Kab)$ for
$v \in W^{k+1,p}(\Kab)$, and obtain the following corollary.
\begin{center}
\fbox{
\begin{minipage}{15truecm}
\begin{corollary}\label{cor23}
For arbitrary $v \in W^{k+1,p}(\Kab)$ $(1 \le p \le \infty)$, the
following estimations hold: 
\begin{align*}
  |v - \I_{\Kab}^k v|_{m,p,\Kab} & \le C_{k,m,p}
   \left(\max\{\alpha,\beta\}\right)^{k+1-m} |v|_{k+1,p,\Kab}.
\end{align*}
\end{corollary}
\end{minipage}
}
\end{center}

The manner of the proof of Theorem~\ref{squeezingtheorem}
is exactly similar as in the previous section.
The ratio $|v|_{m,p,\Kab}^p/|v|_{k+1,p,\Kab}^p$ is written using the
seminorms of $u$ on $\hK$, and is bounded by a constant that does not
depend on $v$.

First, let $1 \le p < \infty$.
For a multi-index $\gamma = (a, b) \in \N_0^2$ and a real $t \neq 0$,
set $(\alpha,\beta)^{\gamma t} := \alpha^{a t}\beta^{b t}$.
Then, we have
{\allowdisplaybreaks
\begin{align}
  \frac{|v|_{m,p,\Kab}^p}{|v|_{k+1,p,\Kab}^p}
   & = \frac{\sum_{|\gamma| = m} \frac{m!}{\gamma!}
      (\alpha,\beta)^{-\gamma p}\left|\partial^{\gamma}u\right|_{0,p,\hK}^p}
      {\sum_{|\delta|= k+1} \frac{(k+1)!}{\delta!}
      (\alpha,\beta)^{- \delta p} \left|\partial^{\delta}
              u \right|_{0,p,\hK}^p} \notag \\
  & = \frac{\sum_{|\gamma| = m}\frac{m!}{\gamma!}
             (\alpha,\beta)^{-\gamma p}
         \left|\partial^{\gamma}u \right|_{0,p,\hK}^p }
       {\sum_{|\gamma| = m} \frac{m!}{\gamma!}
        (\alpha,\beta)^{-\gamma p}
       \left(\sum_{|\eta| = k+1-m} \frac{(k+1-m)!}
             {\eta!(\alpha,\beta)^{\eta p} }\left|\partial^{\eta}
           (\partial^{\gamma}u)\right|_{0,p,\hK}^p\right)} \notag \\
  & \le \frac{\left(\max\{\alpha,\beta\}\right)^{(k+1-m)p}
         \sum_{|\gamma| = m}\frac{m!}{\gamma!}
            (\alpha,\beta)^{-\gamma p}
         \left|\partial^{\gamma}u \right|_{0,p,\hK}^p }
       {\sum_{|\gamma| = m} \frac{m!}{\gamma!}
          (\alpha,\beta)^{-\gamma p}
       \left(\sum_{|\eta| = k+1-m} \frac{(k+1-m)!}
             {\eta!}\left|\partial^{\eta}
        (\partial^{\gamma}u)\right|_{0,p,\hK}^p\right)} \notag \\
   & = \left(\max\{\alpha,\beta\}\right)^{(k+1-m)p}
    \frac{\sum_{|\gamma| = m} \frac{m!}{\gamma!}
         (\alpha,\beta)^{- \gamma p}
         \left|\partial^{\gamma}u \right|_{0,p,\hK}^p}
       {\sum_{|\gamma| = m} \frac{m!}{\gamma!}
      (\alpha,\beta)^{-\gamma p}
        \left|\partial^{\gamma}u \right|_{k+1-m,p,\hK}^p}. \label{eq47}
\end{align}
}
Here, we used the fact that, for a multi-index $\eta$, 
$(\alpha,\beta)^{\eta p} \le \left(\max\{\alpha,\beta\}\right)^{|\eta|p}$
and, for a multi-index $\delta$ with $|\delta| = k+1$,
\[
   \frac{|\delta|!}{\delta!} =
   \sum_{\substack{\gamma+\eta = \delta \\ |\gamma|=m, |\eta|=k+1-m}}
   \frac{|\gamma|!}{\gamma!} \frac{|\eta|!}{\eta!}.
\]

For example, if $k=2$, then we see
{\allowdisplaybreaks
\begin{align*}
  \frac{|v|_{1,p,K_\alpha}^p}{|v|_{3,p,K_\alpha}^p}
   & = \frac{\frac{1}{\alpha^{p}}|u_x|_{0}^p
           + \frac{1}{\beta^{p}}|u_y|_{0}^p}
       {\frac{1}{\alpha^{3p}}|u_{xxx}|_{0}^p 
        + \frac{3}{\alpha^{2p}\beta^{p}} |u_{xxy}|_{0}^p +
        \frac{3}{\alpha^{p}\beta^{2p}} |u_{xyy}|_{0}^p +
        \frac{1}{\beta^{3p}} |u_{yyy}|_{0}^p}  \\
    & = \frac{\frac{1}{\alpha^{p}} |u_x|_{p}^p + \frac{1}{\beta^{p}}|u_y|_{p}^p}
       {\frac{1}{\alpha^p}\left(
    \frac{|u_{xxx}|_{p}^p}{\alpha^{2p}} + 
    2 \frac{|u_{xxy}|_{p}^p}{\alpha^{p}\beta^{p}}  + 
     \frac{|u_{xyy}|_{p}^p}{\beta^{2p}}  \right) +
      \frac{1}{\beta^p}\left(
    \frac{|u_{xxy}|_{p}^p}{\alpha^{2p}} + 
     2\frac{|u_{xyy}|_{p}^p}{\alpha^{p}\beta^{p}}  + 
     \frac{|u_{yyy}|_{p}^p}{\beta^{2p}}  \right)} \\
   & \le \frac{ \max\{\alpha^{2p},\beta^{2p}\}
    \left(\frac{1}{\alpha^{p}}|u_x|_{0}^p
             + \frac{1}{\beta^p}|u_y|_{0}^p \right)}
       {\frac{1}{\alpha^{p}}\left(|u_{xxx}|_{0}^p 
       + 2 |u_{xxy}|_{0}^p + |u_{xyy}|_{0}^p\right) 
       + \frac{1}{\beta^{p}} \left( |u_{xxy}|_{0}^p
       + 2 |u_{xyy}|_{0}^p + u_{yyy}|_{0}^p \right)} \\
   & = \max\{\alpha^{2p},\beta^{2p}\} 
    \frac{\frac{1}{\alpha^{p}}|u_x|_{0,p,\hK}^p
    + \frac{1}{\beta^p}|u_y|_{0,p,\hK}^p}
       {\frac{1}{\alpha^{p}}|u_{x}|_{2,p,\hK}^p
    + \frac{1}{\beta^{p}}|u_{y}|_{2,p,\hK}^p },
\end{align*}
}
and
{\allowdisplaybreaks
\begin{align*}
  \frac{|v|_{2,p,\Kab}^p}{|v|_{3,p,\Kab}^p}
      & = \frac{\frac{1}{\alpha^{2p}}|u_{xx}|_{0}^p
          + \frac{2}{\alpha^p\beta^p}|u_{xy}|_{0}^p
             + \frac{1}{\beta^{2p}}|u_{yy}|_{0}^p}
       {\frac{1}{\alpha^{3p}} |u_{xxx}|_{0}^p
      + \frac{3}{\alpha^{2p}\beta^p} |u_{xxy}|_{0}^p +
        \frac{3}{\alpha^{p}\beta^{2p}} |u_{xyy}|_{0}^p +
       \frac{1}{\beta^{3p}} |u_{yyy}|_{0}^p}  \\
   & = \frac{ \frac{1}{\alpha^{2p}}|u_{xx}|_{p}^p
            + \frac{2}{\alpha^p\beta^p}|u_{xy}|_{p}^p
             + \frac{1}{\beta^{2p}}|u_{yy}|_{p}^p}
       {\frac{1}{\alpha^{2p}}
      \left( \frac{|u_{xxx}|_{p}^p}{\alpha^p} +
            \frac{|u_{xxy}|_{p}^p}{\beta^{p}} \right) +
      \frac{2}{\alpha^p\beta^p}
      \left( \frac{|u_{xxy}|_{p}^p}{\alpha^p} +
            \frac{|u_{xyy}|_{p}^p}{\beta^{p}} \right) +
        \frac{1}{\beta^{2p}} 
               \left( \frac{|u_{xyy}|_{p}^p}{\alpha^p} +
            \frac{|u_{yyy}|_{p}^p}{\beta^{p}} \right)} \\
   & \le \frac{\max\{\alpha^p,\beta^p\}\left(
             \frac{1}{\alpha^{2p}} |u_{xx}|_{0}^p
             + \frac{2}{\alpha^p\beta^p}|u_{xy}|_{0}^p
             + \frac{1}{\beta^{2p}}|u_{yy}|_{0}^p\right)}
       {\frac{1}{\alpha^{2p}} 
     \left(|u_{xxx}|_{0}^p + |u_{xxy}|_{0}^p\right) +
    \frac{2}{\alpha^{p}\beta^p} \left(|u_{xxy}|_{0}^p
       + |u_{xyy}|_{0}^p\right) + \frac{1}{\beta^{2p}}
   \left(|u_{xyy}|_{0}^p + |u_{yyy}|_{0}^p \right)} \\
   & = \max\{\alpha^p,\beta^p\}
    \frac{ \frac{1}{\alpha^{2p}} |u_{xx}|_{0,p,\hK}^p
          + \frac{2}{\alpha^p\beta^p}|u_{xy}|_{0,p,\hK}^p
             + \frac{1}{\beta^{2p}}|u_{yy}|_{0,p,\hK}^p}
       { \frac{1}{\alpha^{2p}} |u_{xx}|_{1,p,\hK}^p +
    \frac{2}{\alpha^{p}\beta^p} |u_{xy}|_{1,p,\hK}^p
       + \frac{1}{\beta^{2p}} |u_{yy}|_{1,p,\hK}^p}.
\end{align*}
}
In the above, we use the notation $|\cdot|_{0}$ instead of
$|\cdot|_{0,p,\hK}$ for simplicity. 

\vvskip
\noindent
\textbf{Exercise:} Confirm the details of the above inequalities,
in particular, \eqref{eq47}.

\vvskip
Now suppose that, for $\T_p^k(\hK)$ and a multi-index $\gamma$,
the set $\Xi_p^{\gamma,k}$ is defined so that
\begin{equation}
    u \in \T_p^k(\hK) \Longrightarrow \partial^{\gamma}
       u \in \Xi_p^{\gamma,k}
   \label{cond1}
\end{equation}
and 
\begin{equation}
   A_p^{\gamma,k} := \sup_{v \in \Xi_p^{\gamma,k}}\frac{|v|_{0,p,\hK}}
   {|v|_{k+1-|\gamma|,p,\hK}} < \infty
   \label{est1}
\end{equation}
hold.  Then, from \eqref{eq47}, we would conclude that
{\allowdisplaybreaks
\begin{align}
  \frac{|v|_{m,p,\Kab}^p}{|v|_{k+1,p,\Kab}^p}
   & \le \left(\max\{\alpha,\beta\}\right)^{(k+1-m)p} 
   \frac{\sum_{|\gamma| = m} \frac{m!}{\gamma!}(\alpha,\beta)^{-\gamma p}
         \left|\partial^{\gamma}u \right|_{0,p,\hK}^p}
       {\sum_{|\gamma| = m} \frac{m!}{\gamma!}
      (\alpha,\beta)^{-\gamma p}
        \left|\partial^{\gamma}u \right|_{k+1-m,p,\hK}^p} \notag\\
   & \le \left(\max\{\alpha,\beta\}\right)^{(k+1-m)p}
   \frac{ \sum_{|\gamma| = m} \left(A_p^{\gamma,k}\right)^{p}
      \frac{m!}{\gamma!}(\alpha,\beta)^{- \gamma p}
         \left|\partial^{\gamma}u \right|_{k+1-m,p,\hK}^p}
       {\sum_{|\gamma| = m} \frac{m!}{\gamma!}
      (\alpha,\beta)^{-\gamma p}
   \left|\partial^{\gamma}u \right|_{k+1-m,p,\hK}^p} \notag \\
%
%
& \le \left(\max\{\alpha,\beta\}\right)^{(k+1-m)p} C_{k,m,p}^p,
 \qquad C_{k,m,p} := \max_{|\gamma|=m} A_p^{\gamma,k}.
  \label{eq30}
\end{align}
}
Our task now is to define $\Xi_p^{\gamma,k}$ that satisfies
\eqref{cond1} and \eqref{est1}.  We will explain the details in the following
sections.

\section{Difference quotients}\label{DQ}
In this section, we define the difference quotients for two-variable
functions.  Our treatment is based on the theory of difference quotients
of one-variable functions given in standard textbooks such as
\cite{Atkinson} and \cite{Yamamoto1}.  All statements in this section
can be readily proved.

\subsection{Difference quotients of one-variable functions}
For a function $f(x)$ and nodal points $x_0, x_1, \cdots, x_n\in\R$,
the difference quotients of $f$ are defined recursively by
\begin{align*}
  & f[x_0,x_1]  := \frac{f(x_0) - f(x_1)}{x_0 - x_1}, \qquad
   f[x_0,x_1,x_2] := \frac{f[x_0,x_1] - f[x_1,x_2]}{x_0 - x_2}, \\
  & f[x_0,x_1,\cdots,x_m] := \frac{f[x_0,\cdots,x_{m-1}]
   - f[x_1,\cdots,x_m]}{x_0 - x_m}.
\end{align*}

A simplest case is $x_i := x_0 + h i$, $i = 1, \cdots, m$, with $h > 0$.
In this case, the difference quotients are 
\begin{align*}
  f[x_0,x_1]  := \frac{f(x_1) - f(x_0)}{h}, \qquad
   f[x_0,x_1,x_2] := \frac{f(x_0) - 2f(x_1) + f(x_2)}{h^2},
\end{align*}
and so on.  The difference quotients are expressed by integration:
\begin{align*}
 f[x_0,x_1] & = \frac{1}{x_1 - x_0}\int_{x_0}^{x_1} f'(t) \dd t
   = \int_0^1 f'(x_0 + t_1 (x_1 - x_0)) \dd t_1, \\
 f[x_0,x_1,x_2] & 
 = \frac{f[x_2,x_1] - f[x_1,x_0]}{x_2 - x_0}
 = \frac{f[x_2,x_0] - f[x_0,x_1]}{x_2 - x_1} \\
 & = \frac{1}{x_2-x_1}\int_0^1   \left( f'(x_0 + t_1(x_2 - x_0))
   - f'(x_0 + t_1(x_1 - x_0)) \right)\dd t_1 \\
  & = \int_0^1\int_0^{t_1}
  f''\left(x_0 + t_1(x_1-x_0) + t_2(x_2-x_1)\right) \dd t_2\dd t_1.
\end{align*}
For $n \ge 1$, the following formula holds:
\begin{align}
 f[x_0,x_1,\cdots,x_n] & = \int_0^1\int_0^{t_1}\cdots\int_0^{t_{n-1}}
  f^{(n)}\left(x_0 + \sum_{i=1}^n t_i(x_i-x_{i-1})\right)
   \dd t_n \cdots \dd t_2 \dd t_1.
   \label{int-quo}
\end{align}

\noindent
\textbf{Exercise:} Prove \eqref{int-quo} by induction.

\subsection{Difference quotients of two variable functions}
We now extend the difference quotient to functions with two
variables. For a positive integer $k$, the set
$\widehat{\Sigma}^k \subset \hK$ is defined by
\begin{align*}
   \widehat{\Sigma}^k & := \Sigma^k(\hK) :=\left\{ \bfx_{\gamma}
   := \frac{\gamma}{k}
     \in \hK \biggm|  \gamma \in \N_0^2, \; 0 \le |\gamma| \le k \right\},
\end{align*}
where $\gamma/k=(a_1/k,a_2/k)$ is understood as the
coordinate of a point in $\widehat{\Sigma}^k$.

For $\bfx_{\gamma} \in \widehat{\Sigma}^k$ and a multi-index $\delta \in \N_0^2$
with $|\gamma| \le k - |\delta|$, we define the correspondence $\Delta^\delta$
between nodes by
\[
   \Delta^\delta \bfx_{\gamma} := \bfx_{\gamma+\delta} \in \widehat{\Sigma}^k.
\]
For example, $\Delta^{(1,1)}\bfx_{(0,0)} = \bfx_{(1,1)}$ and
$\Delta^{(2,1)}\bfx_{(0,1)} = \bfx_{(2,2)}$.
Using $\Delta^\delta$, we define the \textbf{difference quotients}
on $\widehat{\Sigma}^k$ for $f \in C^0(\hK)$ by
\begin{align*}
 f^{|\delta|}[\bfx_{\gamma},\Delta^\delta \bfx_{\gamma}] :=
  k^{|\delta|}\sum_{\eta \le \delta} 
\frac{(-1)^{|\delta|-|\eta|}}{\eta!(\delta-\eta)!}
    f(\Delta^\eta \bfx_{\gamma}).
\end{align*}
For simplicity, we denote 
$f^{|\delta|}[\bfx_{(0,0)}, \Delta^\delta\bfx_{(0,0)}]$ by
$f^{|\delta|}[\Delta^\delta\bfx_{(0,0)}]$.  The following are
examples of $f^{|\delta|}[\Delta^\delta\bfx_{(0,0)}]$: 
\begin{align*}
   f^2[\Delta^{(2,0)}\bfx_{(0,0)}] & =
   \frac{k^2}{2} (f(\bfx_{(2,0)}) - 2 f(\bfx_{(1,0)})
     + f(\bfx_{(0,0)})), \\
   f^2[\Delta^{(1,1)}\bfx_{(0,0)}] & =
   k^2 (f(\bfx_{(1,1)}) -  f(\bfx_{(1,0)})
    - f(\bfx_{(0,1)}) + f(\bfx_{(0,0)})), \\
  f^3[\Delta^{(2,1)}\bfx_{(0,0)}] & =
   \frac{k^3}{2} (f(\bfx_{(2,1)}) - 2 f(\bfx_{(1,1)})
   + f(\bfx_{(0,1)}) - f(\bfx_{(2,0)}) \\
 & 
  \hspace{2.42truecm}
   + 2 f(\bfx_{(1,0)}) - f(\bfx_{(0,0)})).
\end{align*}
Let $\eta \in \N_0^2$ be such that $|\eta|=1$ and $\eta \le \delta$.
The difference quotients clearly satisfy the following
recursive relations:
\begin{align*}
 f^{|\delta|}[\bfx_{\gamma},\Delta^{\delta}\bfx_{\gamma}] & = 
 \frac{k}{\delta\cdot\eta}  \left(
 f^{|\delta|-1}[\bfx_{\gamma+\eta},\Delta^{\delta-\eta}\bfx_{\gamma+\eta}]
-  f^{|\delta|-1}[\bfx_{\gamma},\Delta^{\delta-\eta}\bfx_{\gamma}]  \right).
\end{align*}

If $f \in C^k(\hK)$, the difference quotient
$f^{|\delta|}[\bfx_{\gamma},\Delta^\delta \bfx_{\gamma}]$ is written as
an integral of $f$.  Setting $d=2$ and $\delta=(0,s)$, for example,
we have
\begin{gather*}
   f^{1}[\bfx_{(l,q)},\Delta^{(0,1)} \bfx_{(l,q)}]
   = k(f(\bfx_{(l,q+1)}) - f(\bfx_{(l,q)}))
   = \int_{0}^{1} \partial^{(0,1)}
    f\left(\frac{l}{k},\frac{q}{k} + \frac{w_1}{k}\right) \dd w_1, \\
  \hspace{-3truecm}
  f^{1}[\bfx_{(l,q)},\Delta^{(0,2)} \bfx_{(l,q)}]
   = \frac{k^2}{2}(f(\bfx_{(l,q+2)}) - 2 f(\bfx_{(l,q+1)}) + f(\bfx_{(l,q)})) \\
  \hspace{1.9truecm}
 = \int_{0}^{1}\int_0^{w_1} \partial^{(0,2)}
  f\left(\frac{l}{k},\frac{q}{k} + \frac{1}{k}(w_1 + w_2)
   \right) \dd w_2 \dd w_1 \\
  \hspace{4.2truecm}
  = k \int_{0}^{1} \left[
   \partial^{(0,1)}  f\left(\frac{l}{k},\frac{q}{k} + \frac{2}{k}w_1 \right) 
  - \partial^{(0,1)} f\left(\frac{l}{k},\frac{q}{k} + \frac{1}{k}w_1\right)
    \right]\dd w_1, \\
 \hspace{-10.5truecm}
f^{s}[\bfx_{(l,p)},\Delta^{(0,s)} \bfx_{(l,q)}]  \\
\hspace{1truecm}
= \int_{0}^{1}\int_0^{w_1}\cdots \int_0^{w_{s-1}} \partial^{(0,s)}
  f\left(\frac{l}{k},\frac{q}{k} + \frac{1}{k}(w_1 + \cdots + w_s)
   \right) \dd w_s \cdots \dd w_2 \dd w_1.
\end{gather*}

To provide a concise expression for the above integral, we introduce
the $s$-simplex
\begin{gather*}
   \Simp_s := \left\{(t_1,t_2,\cdots,t_s) \in \R^s \mid 
   t_i \ge 0, \ 0 \le t_1 + \cdots + t_s \le 1 \right\},
\end{gather*}
and the integral of $g \in L^{1}(\Simp_s)$ on $\Simp_s$ is defined by
\begin{gather*}
  \int_{\Simp_s} g(w_1,\cdots,w_k) \dd\mathbf{W_s}
  := \int_{0}^{1}\int_0^{w_1}\cdots \int_0^{w_{s-1}} 
  g(w_1, \cdots, w_s)  \dd w_s \cdots \dd w_2 \dd w_1,
\end{gather*}
where $\dd \mathbf{W_s} = \dd w_s \cdots \dd w_2\dd w_1$.
Then, $f^{s}[\bfx_{(l,q)},\Delta^{(0,s)} \bfx_{(l,q)}]$ becomes
\begin{align*}
f^{s}[\bfx_{(l,q)},\Delta^{(0,s)} \bfx_{(l,q)}]
  & = \int_{\Simp_s} \partial^{(0,s)}
  f\left(\frac{l}{k}, \mathbf{W_s}\right) \dd \mathbf{W_s}, \quad
   \mathbf{W_s} := \frac{q}{k} + \frac{1}{k}(w_1 + \cdots + w_s).
\end{align*}
For a general multi-index $(t,s)$, we have
\begin{gather*}
 f^{t+s}[\bfx_{(l,q)},\Delta^{(t,s)} \bfx_{(l,q)}]
   = \int_{\Simp_s}\int_{\Simp_t}\partial^{(t,s)}
  f\left(\mathbf{Z_t},\mathbf{W_s}\right) \dd \mathbf{Z_t}
    \dd \mathbf{W_s}. \\
  \mathbf{Z_t} := \frac{l}{k} + \frac{1}{k}(z_1 + \cdots + z_t), \quad
  \dd\mathbf{Z_t} := \dd z_t \cdots \dd z_2\dd z_1.
\end{gather*}

Let $\square_{\gamma}^\delta$ be the rectangle defined by 
$\bfx_{\gamma}$ and $\Delta^\delta \bfx_{\gamma}$
as the diagonal points.  If $\delta=(t,0)$
or $(0,s)$, $\square_{\gamma}^\delta$ degenerates to a segment.  For
$v \in W^{1,1}(\hK)$ and $\square_{\gamma}^\delta$ with $\gamma=(l,q)$, we
denote the integral as
\begin{equation*}
   \int_{\square_{\gamma}^{(t,s)}} v :=
   \int_{\Simp_s}\int_{\Simp_t}
     v\left(\mathbf{Z_t},\mathbf{W_s}\right)
    \dd \mathbf{Z_t} \dd \mathbf{W_s}.
\end{equation*}
If $\square_{\gamma}^\delta$ degenerates to a segment, the integral is
understood as an integral on the segment.  By this notation, the
difference quotient $f^{t+s}[\bfx_{\gamma},\Delta^{(t,s)} \bfx_{\gamma}]$ is
written as 
\begin{align*}
 f^{(t+s)}[\bfx_{\gamma},\Delta^{(t,s)} \bfx_{\gamma}] 
  = \int_{\square_{\gamma}^{(t,s)}} \partial^{(t,s)} f.
\end{align*}
Therefore, if $u \in \T_p^k(\hK)$, then we have
\begin{align}
 0 = u^{t+s}[\bfx_{\gamma},\Delta^{(t,s)} \bfx_{\gamma}] 
  = \int_{\square_{\gamma}^{(t,s)}} \partial^{(t,s)} u, \qquad
  \forall \square_{\gamma}^{(t,s)} \subset \hK.
  \label{tomoko}
\end{align}

\vvskip
\noindent
\textbf{Exercise: } Confirm that all the equations in this section
certainly hold.

\section{The proof of Theorem~\ref{squeezingtheorem}}
By introducing the notation in the previous section, we now be able to
define $\Xi_p^{\gamma,k} \subset W^{k+1-|\gamma|,p}(\hK)$ and $A_p^{\gamma,k}$
for $p \in [1,\infty]$, which satisfy \eqref{cond1} and \eqref{est1}.
For multi-index $\gamma$, define
\begin{align*}
  \Xi_p^{\gamma,k} & := \left\{ v \in W^{k+1-|\gamma|,p}(\hK) \Bigm| 
    \int_{\square_{lp}^{\gamma}} v = 0,\quad \forall
   \square_{lp}^{\gamma} \subset \hK \right\}.
\end{align*}
From the definition and \eqref{tomoko}, it is clear that \eqref{cond1} holds.
Define 
\begin{align*}
   A_p^{\gamma,k} := \sup_{v \in \Xi_p^{\gamma,k}} \frac{|v|_{0,p,\hK}}
  {|v|_{k+1-|\gamma|,p,\hK}}, \qquad 1 \le p \le \infty.
\end{align*}
Then, the following lemma holds.

\begin{center}
\fbox{
\begin{minipage}{15truecm}
\begin{lemma}\label{chiaki} We have
$\Xi_p^{\gamma,k} \cap \PP_{k-|\gamma|} = \{0\}$.  That is,
if $q \in \PP_{k-|\gamma|}$ belongs to $\Xi_p^{\gamma,k}$, then $q=0$.
\end{lemma}
\end{minipage}
}
\end{center}
\textbf{Proof:}
 We notice that
$\mathrm{dim} \PP_{k-|\delta|}
= \#\{\square_{lp}^\delta \subset \hK \}$.
For example, if $k=4$ and $|\delta|=2$, then
$\mathrm{dim}\PP_2 = 6$.
This corresponds to the fact that,
in $\hK$, there are six squares with size $1/4$  for $\delta=(1,1)$ and
there are six horizontal segments of length $1/2$ for $\delta=(2,0)$.
All their vertices (corners and end-points) belong to $\Sigma^4(\hK)$
(see Figure~\ref{box_relation}).  Now, suppose that $v \in \PP_{k-|\delta|}$
satisfies $\int_{\square_{lp}^{\delta}} q = 0$ for all
$\square_{lp}^{\delta} \subset \hK$.  This condition is linearly
independent and determines $q = 0$ uniquely. $\square$

\newcommand{\miho}{1.1pt}
\vspace{3mm}
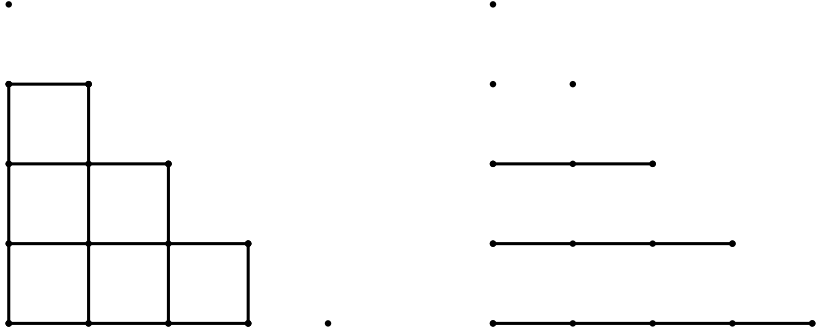
\begin{figure}[thb]
\begin{center}
\begin{tikzpicture}[line width = 1pt,scale=2]
\coordinate (P_1) at (0.0,0.0);
\coordinate (P_2) at (2.0,0.0);
\coordinate (P_3) at (0.0,2.0);
\coordinate (P_4) at (1.5,0.0);
\coordinate (P_5) at (1.5,0.5);
\coordinate (P_6) at (1.0,0.0);
\coordinate (P_7) at (1.0,0.5);
\coordinate (P_8) at (1.0,1.0);
\coordinate (P_9) at (0.5,0.0);
\coordinate (P_10) at (0.5,0.5);
\coordinate (P_11) at (0.5,1.0);
\coordinate (P_12) at (0.5,1.5);
\coordinate (P_13) at (0.0,0.5);
\coordinate (P_14) at (0.0,1.0);
\coordinate (P_15) at (0.0,1.5);
\draw (P_15)--(P_1)--(P_4);
\draw (P_15)--(P_12)--(P_9);
\draw (P_14)--(P_8)--(P_6);
\draw (P_13)--(P_5)--(P_4);
\fill (P_1) circle (\miho) (P_2) circle (\miho) (P_3) circle (\miho)
 (P_4) circle (\miho) (P_5) circle (\miho) (P_6) circle (\miho) (P_7)
 circle (\miho) (P_8) circle (\miho) (P_9) circle (\miho) (P_10) circle
 (\miho) (P_11) circle (\miho) (P_12) circle (\miho) (P_13) circle
 (\miho) (P_14) circle (\miho) (P_15) circle (\miho);
\end{tikzpicture}
  \hspace{1cm}
\qquad
\begin{tikzpicture}[line width = 1pt,scale=2]
\coordinate (P_1) at (0.0,0.0);
\coordinate (P_2) at (2.0,0.0);
\coordinate (P_3) at (0.0,2.0);
\coordinate (P_4) at (1.5,0.0);
\coordinate (P_5) at (1.5,0.5);
\coordinate (P_6) at (1.0,0.0);
\coordinate (P_7) at (1.0,0.5);
\coordinate (P_8) at (1.0,1.0);
\coordinate (P_9) at (0.5,0.0);
\coordinate (P_10) at (0.5,0.5);
\coordinate (P_11) at (0.5,1.0);
\coordinate (P_12) at (0.5,1.5);
\coordinate (P_13) at (0.0,0.5);
\coordinate (P_14) at (0.0,1.0);
\coordinate (P_15) at (0.0,1.5);
\coordinate (a) at (0.0,0.489);
\coordinate (b) at (1.0,0.489);
\coordinate (c) at (0.5,0.511);
\coordinate (d) at (1.5,0.511);
\coordinate (e) at (0.0,-0.011);
\coordinate (f) at (2.0,-0.011);
\coordinate (g) at (0.5,0.011);
\coordinate (h) at (1.5,0.011);
\draw (P_14)--(P_8);
\draw (a)--(b);
\draw (c)--(d);
\draw (e)--(f);
\draw (g)--(h);
\fill (P_1) circle (\miho) (P_2) circle (\miho) (P_3) circle (\miho)
 (P_4) circle (\miho) (P_5) circle (\miho) (P_6) circle (\miho) (P_7)
 circle (\miho) (P_8) circle (\miho) (P_9) circle (\miho) (P_10) circle
 (\miho) (P_11) circle (\miho) (P_12) circle (\miho) (P_13) circle
 (\miho) (P_14) circle (\miho) (P_15) circle (\miho);
\end{tikzpicture}
\caption{The six squares of size $1/4$ for $\delta=(1,1)$
and the (union of) six segments of length $1/2$ for
$\delta=(2,0)$ in $\hK$.}
 \label{box_relation}
\end{center}
\end{figure}

\vvskip
To understand the above proof clearly, we consider the cases
$k=2$ and $3$.   Let $k=2$ and $\gamma=(1,0)$.  Then, 
$k-|\gamma|=1$.   Set $q(x,y)=a+bx+cy$.  If the
three integrals
\begin{align*}
  \int_{\square_{00}^{(1,0)}} q(x,y) = a + \frac{b}{4},
  \qquad 
  \int_{\square_{10}^{(1,0)}} q(x,y) = a + \frac{3b}{4},
  \qquad 
  \int_{\square_{01}^{(1,0)}} q(x,y) = a + \frac{b}{4} + \frac{c}{2}
\end{align*}
are equal to $0$, then we have $a=b=c=0$, that is, $q(x,y)=0$.
The case $\gamma=(0,1)$ is similar.

Let $k=3$ and $\gamma=(1,0)$.  Then,
$k-|\gamma|=2$. Set $q(x,y)=a+bx+cy+dx^2+ey^2+fxy$.
If the integrals
\begin{align*}
  \int_{\square_{00}^{(1,0)}} q(x,y) & = a + \frac{b}{6} + \frac{d}{27},
  \qquad 
  \int_{\square_{10}^{(1,0)}} q(x,y) = a + \frac{b}{2} + \frac{7}{27}d, \\
  \qquad 
  \int_{\square_{20}^{(1,0)}} q(x,y) & = a + \frac{5}{6}b + \frac{19}{27}d
\end{align*}
are all equal to $0$, we have $a=b=d=0$.  Moreover, if the integrals
{\allowdisplaybreaks
\begin{align*}
  \int_{\square_{01}^{(1,0)}} q(x,y) & = \frac{c}{3} + \frac{e}{9} + \frac{f}{18},
  \qquad 
  \int_{\square_{11}^{(1,0)}} q(x,y) = \frac{c}{3} + \frac{e}{9} + \frac{f}{6}, \\
  \qquad 
  \int_{\square_{02}^{(1,0)}} q(x,y) & = \frac{2}{3}c + \frac{4}{9}e + \frac{f}{9}
\end{align*}
}
are equal to $0$ as well, we have $c=e=f=0$.  Hence, we conclude that
$q(x,y)=0$.  The case $\gamma=(0,1)$ is similar.


\vvskip
\begin{center}
\fbox{
\begin{minipage}{15truecm}
\begin{lemma}\label{lem22} We have $A_p^{\gamma,k} < \infty$,
 $p \in [1,\infty]$.  That is, \eqref{est1} holds.
\end{lemma}
\end{minipage}
}
\end{center}
\textbf{Proof:} The proof is by contradiction.
Suppose that $A_p^{\gamma,k} = \infty$.  Then, there exists a sequence 
$\{u_n\}_{n=1}^\infty \subset \Xi_p^{\gamma,k}$ such that
\[
   |u_n|_{0,p,\hK} = 1, \qquad
   \lim_{n \to \infty} |u_n|_{k+1-|\gamma|,p,\hK} = 0.
\]
By the inequality \eqref{eq3.1.9}, for an arbitrary $\varepsilon > 0$,
there exists a sequence $\{q_n\}\subset \PP_{k-|\gamma|}$ such that
\begin{align*}
   \inf_{q \in \PP_{k-|\gamma|}}
  \|u_n + q\|_{k+1-|\gamma|,p,\hK} 
  & \le \|u_n + q_n\|_{k+1-|\gamma|,p,\hK} \\
  & \le \inf_{q \in \PP_{k-|\gamma|}}
  \|u_n + q\|_{k+1-|\gamma|,p,\hK} + \frac{\varepsilon}{n} \\
  & \le C |u_n|_{k+1-|\gamma|,p,\hK} + \frac{\varepsilon}{n}, \quad
   \lim_{n \to \infty} \|u_n + q_n\|_{k+1-|\gamma|,p,\hK} = 0.
\end{align*}
Since $|u_n|_{k+1-|\gamma|,p,\hK}$ and $|u_n|_{0,p,\hK}$ are bounded,
$|u_n|_{m,p,\hK}$ ($1 \le m \le k-|\gamma|$) is bounded as well by
Gagliardo--Nirenberg's inequality (Theorem~\ref{GagliardoNirenberg}).
That is, $\|u_n\|_{k+1-|\gamma|,p,\hK}$ and
$\{q_n\} \subset \PP_{k-|\gamma|}$ are bounded.  Thus, there exists
a subsequence $\{q_{n_i}\}$ such that $q_{n_i}$ converges to
$\bar{q} \in \PP_{k-|\gamma|}$.  In particular, we see
\[
   \lim_{n_i \to \infty} \|u_{n_i} + \bar{q}\|_{k+1-|\gamma|,p,\hK} = 0.
\]
Therefore, for any $\square_{lp}^\gamma$, we notice that
\[
  0 = \lim_{n_i\to\infty} \int_{\square_{lp}^\gamma}
   (u_{n_i} + \bar{q}) = \int_{\square_{lp}^\gamma} \bar{q},
\]
and $\bar{q} = 0$ by Lemma~\ref{chiaki}.  This yields
\[
  \lim_{n_i \to \infty}\|u_{n_i}\|_{k+1-|\gamma|,p,\hK} = 0,
\]
which contradicts 
$\lim_{n_i \to \infty}\|u_{n_i}\|_{k+1-|\gamma|,p,\hK} \ge
 \lim_{n_i \to \infty}|u_{n_i}|_{0,p,\hK}=1$. $\square$

\vvskip
Now, we have defined the set $\Xi_p^{\gamma,k}$ that satisfies
\eqref{cond1} and the estimate \eqref{est1} has been shown. 
Therefore, 
Theorem~\ref{squeezingtheorem} has been proved by \eqref{eq30}.

\vvskip
\noindent
\textbf{Exercise: } We have shown the Theorem~\ref{squeezingtheorem}
for the case $1 \le p < \infty$.  Prove Theorem~\ref{squeezingtheorem}
for the case $p=\infty$.

\section{The error estimation on general triangles in terms of circumradius}
Using the previous results, we can obtain the error estimations
on general triangles.  Recall the reference triangle and the definition
of the standard position of an aribtrary triangle $K$
(Figure~\ref{LiuKikuchiTriangle0}).
Let $\Kab$ be the triangle with the vertices
$(0,0)^{\top}$, $(\alpha,0)^{\top}$, and $(0,\beta)^{\top}$.  Let
$\hK$ be the reference triangle with the vertices 
$(0,0)^{\top}$, $(1,0)^{\top}$, and $(0,1)^{\top}$.

\begin{figure}[thb]
\begin{center}
\begin{tikzpicture}[line width = 1pt,scale=0.6]
   \coordinate [label=below:{$\bfx_1$}](A) at (0.0,0.0);
   \coordinate [label=below:{$\bfx_2$}](B) at (7.0,0.0);
   \coordinate [label=above:{$\bfx_3$}](C) at (-2.0,4.0);
   \draw (A) -- node[below]{$\alpha$} (B) ;
   \draw (B) -- node[pos=0.4,above]{$h_K$} (C) ;
   \draw (C) -- node[pos=0.6, left]{$\beta$} (A) ;
   \coordinate (D) at (0.8,0.0);
   \coordinate (E) at ($(A)!0.17!(C)$);
   \draw [bend right,thin] (D) to node[pos=0.3,above]{$\theta$} (E) ;
   \coordinate [label=:{$K$}](F) at (2.0,0.57);
\end{tikzpicture}
 \caption{The standard position of a general triangle (reprint).
 The vertices are
 $\bfx_1=(0,0)^\top$, $\bfx_2=(\alpha,0)^\top$, and
  $\bfx_3=(\beta s,\beta t)^\top$,
 where $s^2 + t^2 = 1$, $t > 0$.
 We assume that $0 < \beta \le \alpha \le h_K$.  Then,
$\pi/3 \le \theta < \pi$. }
 \label{LiuKikuchiTriangle}
\end{center}
\end{figure}
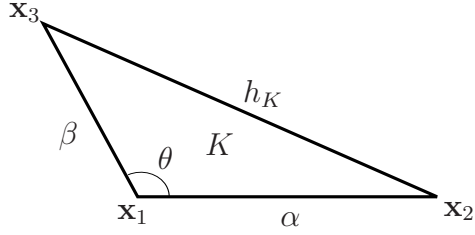
We consider $2\times 2$ matrices
\begin{align*}
    A & := \begin{pmatrix}
	   \alpha & \beta s \\ 0 & \beta t 
	 \end{pmatrix}
      = \begin{pmatrix}
	   1 &  s \\ 0 &  t 
	 \end{pmatrix}
        \begin{pmatrix}
	   \alpha &  0 \\ 0 &  \beta 
	 \end{pmatrix}, \quad
  \wA := \begin{pmatrix}
	   1 &  s \\ 0 &  t 
	 \end{pmatrix}, \quad
   D_{\alpha\beta} := \begin{pmatrix}
	   \alpha &  0 \\ 0 & \beta 
	 \end{pmatrix},    \\
  & A^{-1} = \begin{pmatrix}
	   \alpha^{-1} & -\alpha^{-1}st^{-1} \\ 0 & \beta^{-1}t^{-1}
	 \end{pmatrix}
    = \begin{pmatrix}
	 \alpha^{-1} & 0 \\ 0 & \beta^{-1}
	 \end{pmatrix}
      \begin{pmatrix}
	   1 & - st^{-1} \\ 0 & t^{-1}
	 \end{pmatrix},
\end{align*}
and the linear transformation $\bfy = A\bfx$. The reference triangle
$\hK$ is transformed to $\Kab$ by $\bfy = D_{\alpha\beta}\bfx$, and
$\Kab$ is transformed to $K$ by $\bfy = \wA\bfx$.  Accordingly, 
$\T_p^k(K)$ is pulled-back to $\T_p^k(\Kab)$ by the mapping
$\T_p^k(K) \ni v \mapsto \hv := v\circ \wA$, and 
$\T_p^k(\Kab)$ is pulled-back to $\T_p^k(\hK)$ by the mapping
$\T_p^k(K) \ni v \mapsto \hv := v\circ D_{\alpha\beta}$.

By Theorem~\ref{squeezingtheorem}, for arbitrary $\alpha \ge\beta > 0$
and arbitrary $p$, $1 \le p \le \infty$, there exists a constant
$C_{k,m,p}$ depending only on $k$, $m$, $p$ such that 
\begin{equation}
   B_p^{m,k}(\Kab) := \sup_{v \in \T_p^{k}(\Kab)}
    \frac{|v|_{m,p,\Kab}}{|v|_{k+1,p,\Kab}}
  \le \alpha^{k+1-m} C_{k,m,p}.
  \label{kobayashi-tsuchiya}
\end{equation}

A simple computation confirms that $\wA^\top\wA$ has the eigenvalues 
$1 \pm |s|$, and $\widetilde{A}^{-1}\widetilde{A}^{-\top}$ has
the eigenvalues $(1\pm|s|)^{-1}$.  That is,
$\|\widetilde{A}\| = (1 + |s|)^{1/2}$, 
$\|\widetilde{A}^{-1}\| = (1 - |s|)^{-1/2}$, and
$\det \widetilde{A} = t$.  Therefore, defining
$\hv(\bfx) = v(\widetilde{A}\bfx)$ for $v \in \T_p^k(K)$,
it follows from \eqref{generalmp} that
\begin{gather*}
  |v|_{m,p,K} \le 2^{m\mu(p)}  t^{1/p}
    \|\widetilde{A}^{-1}\|^{m} |\hv|_{m,p,\Kab}, \\
  2^{-(k+1) \mu(p)} t^{1/p} \|\widetilde{A}\|^{-(k+1)}
   |\hv|_{k+1,p,\Kab} \le |v|_{k+1,p,K}.
\end{gather*}
Combining the above inequalities and \eqref{kobayashi-tsuchiya}, we
obtain
\begin{align*}
   \frac{|v|_{m,p,K}}{|v|_{k+1,p,K}} & \le
   c_{k,m,p} \|\widetilde{A}\|^{k+1}
    \|\widetilde{A}^{-1}\|^{m}
   \frac{|\hv|_{m,p,\Kab}}{|\hv|_{k+1,p,\Kab}} \\
  & \le c_{k,m,p}C_{k,m,p} \|\widetilde{A}\|^{k+1}
    \|\widetilde{A}^{-1}\|^{m}\alpha^{k+1-m},
\end{align*}
where $c_{k,m,p} := 2^{(k+1+m)\mu(p)}$.  Hence, we obtain the following
lemma.

\begin{center}
 \fbox{
\begin{minipage}{15truecm}
\begin{lemma}  \label{liu-kikuchi}
For an arbitrary triangle $K$ in the standard position, we have 
\begin{align*}
   B_p^{m,k}(K) & \le c_{k,m,p} 
  \|\widetilde{A}\|^{k+1}\|\widetilde{A}^{-1}\|^{m}
  B_p^{m,k}(\Kab) \\
  & \le c_{k,m,p} C_{k,m.p}
  \|\widetilde{A}\|^{k+1}\|\widetilde{A}^{-1}\|^{m} 
   \alpha^{k+1-m},
\end{align*} 
where $\|\widetilde{A}\| = (1 + |s|)^{1/2}$ and 
$\|\widetilde{A}^{-1}\| = (1 - |s|)^{-1/2}$.
\end{lemma}
\end {minipage}
}
\end{center}

Applying Lemma~\ref{liu-kikuchi} to $v - \I_K^k v \in \T_p^k(K)$, we
have the following corollary.
\begin{center}
 \fbox{
\begin{minipage}{15truecm}
\begin{corollary}  
For an arbitrary triangle $K$ in the standard position, we have 
\begin{align*}
 |v - \I_K^k v|_{m,p,K} \le c_{k,m,p} C_{k,m.p}
  \|\widetilde{A}\|^{k+1}\|\widetilde{A}^{-1}\|^{m} 
   \alpha^{k+1-m}|v|_{k+1,p,K}, \quad
   \forall v \in W^{k+1,p}(K).
\end{align*} 
\end{corollary}
\end {minipage}
}
\end{center}

We would like to obtain upper bounds of 
$\|\widetilde{A}\|$ and $\|\widetilde{A}^{-1}\|$.
From Lemma~\ref{liu-kikuchi}, we obviously have
$\|\widetilde{A}\| \le \sqrt{2}$.  For 
$\|\widetilde{A}^{-1}\|$, we observe that
{\allowdisplaybreaks
\begin{align}
   \|\widetilde{A}^{-1}\| = \frac{1}{(1 - |s|)^{1/2}}
   & = \frac{(1+|s|)^{1/2}}{t}
   \qquad (\because\ s^2 + t^2 = 1) \notag  \\
   & \le \frac{2^{1/2}\alpha\beta h_K}{\alpha\beta h_K t} 
   = \frac{\alpha\beta h_K}{2^{1/2}h_K|K|}
   \qquad \left(\because |K| = \frac{1}{2}\alpha \beta t\right) 
 \notag \\
&   = \frac{2^{3/2}R_K}{h_K} 
   \qquad \left(\because R_K = \frac{\alpha\beta h_K}{4|K|}\right).
         \label{circumradius}
\end{align}
}
Thus, redefining the constant $C_{k,m,p}$, we obtain the following
theorem.
\begin{center}
 \fbox{
\begin{minipage}{15truecm}
\begin{theorem}\label{thm18}
Suppose that a triangle $K$ is in the standard position.
Let $k$, $m$ be integers with $k \ge 1$, $m = 0,\cdots,k$ and
$1 \le p \le \infty$.
Then, the following estimate holds:
\begin{align*}
    B_p^{m,k}(K) := \sup_{v\in \T_p^{k}(K)} 
    \frac{|v|_{m,p,K}}{|v|_{k+1,p,K}} \le C_{k,m,p} \,
     \left(\frac{R_K}{h_K}\right)^m \alpha^{k+1-m},
\end{align*}
where $R_K$ is the circumradius of $K$, and $C_{k,m,p}$ is a 
constant depending only on $k$, $m$, and $p$.
\end{theorem}
\end {minipage}
}
\end{center}

\vspace{0.2truecm}
Now, let $K$ be an arbitrary triangle.
Note that $\alpha \le h_K$ and the Sobolev norms
are affected by rotations if $p \neq 2$ up to an constant
(see \eqref{rotationp}). 
Then, with rewriting the constant, we obtain the following corollary
from Theorem~\ref{thm18}, that is the main theorem
of this survey (reprint of Theorem~\ref{eq:circumradius-est}).
\begin{center}
 \fbox{
\begin{minipage}{15truecm}
\begin{corollary}\label{cor21}
Let $K$ be an arbitrary triangle with circumradius $R_K$.
Let $k$ and $m$ be intergers with $k \ge 1$ and $m=0,\cdots,k$.
Let $p$, $1 \le p \le \infty$.
For the Lagrange interpolation $\I_K^k v$ of degree $k$ on $K$,
the following estimate holds: for any $v \in W^{2,p}(K)$,
\begin{align*}
  & B_p^{m,k}(K) := \sup_{u\in \T_p^{k}(K)} 
   \frac{|u|_{m,p,K}}{|u|_{k+1,p,K}} \le \, C_{k,m,p} 
  \left(\frac{R_K}{h_k}\right)^m h_K^{k+1-m}, \\
  & |v - \I_K^k v|_{m,p,K} \le  C_{k,m,p} 
    \left(\frac{R_K}{h_K}\right)^m
   h_K^{k+1-m} |v|_{k+1,p,K},
\end{align*}
where $C_{k,m,p}$ depends only on $k$, $m$, and $p$.
\end{corollary}
\end {minipage}
}
\end{center}

\vvskip
\textit{Remarks:} (1) Let $\Omega \subset \R^2$ be a bounded polygonal
domain. We compute a numerical solution of the Poisson equation
\begin{align*}
    -\Delta u = f \text{ in } \Omega, \quad
    u = 0 \text{ on } \partial\Omega
\end{align*}
by the conforming piecewise $k$th-order finite element method on 
simplicial elements.  To this end, we construct a triangulation
$\T_h$ of $\Omega$ and consider the piecewise $\PP_k$ continuous
function space $S_h \subset H_0^1(\Omega)$.
The weak form of the Poisson equation is
\begin{align*}
    \int_\Omega \nabla u \cdot \nabla v \dd \bfx  = 
    \int_\Omega f v \dd \bfx, \quad \forall
      v \in H_0^1(\Omega),
\end{align*}
and the finite element solution is defined as the unique solution
$u_h \in S_h$ of
\begin{align*}
    \int_\Omega \nabla u_h \cdot \nabla v_h \dd \bfx  = 
    \int_\Omega f v_h \dd \bfx, \quad \forall v_h \in S_h.
\end{align*}
C\'ea's Lemma implies that the error $|u - u_h|_{1,2,\Omega}$ is
estimated as
\begin{align}
  |u - u_h|_{1,2,\Omega} \le \left(
   \sum_{K \in \T_h} |u - \I_K^k u|_{1,2,K}^2 \right)^{1/2}.
  \label{Cea_lemma}
\end{align}
Combining \eqref{Cea_lemma} and Corollary~\ref{cor21} with $p=2$,
$k \ge 2$, $m=1$, we have
\begin{align*}
  |u - u_h|_{1,2,\Omega} & \le C\left(\sum_{K \in \T_h}
   |u - \I_K^ku|_{k+1,2,K}^2\right)^{1/2} \\
  & \le C\left(\sum_{K \in \T_h}
   (R_K h_K^{k-1})^2 |u|_{k+1,2,K}^2\right)^{1/2} \\
  & \le C \max_{K \in \T_h}(R_Kh_K^{k-1}) |u|_{k+1,2,\Omega}. 
\end{align*}
Therefore, if $\max_{K \in \T_h}(R_Kh_K^{k-1}) \to 0$ as $h \to 0$ and
$u \in H^{k+1}(\Omega)$, the finite element solution $u_h$ converges
to the exact solution $u$ even if there exist many skinny elements
violating the shape regularity condition or the maximum angle condition
in $\T_h$. 

Recall the triangle depicted in Figure~\ref{example1} (right) with
vertices $(0,0)^\top$, $(h,0)^\top$, and $(h^\alpha,h^\beta)^\top$
with $R_K = \mathcal{O}(h^{1 + \alpha - \beta})$.
Suppose now that $\alpha + 1 \le \beta < 2 + \alpha$.
If a sequence of triangulations contains those triangles, and
$k = 1$, then $\max_{K \in \T_h} R_K = \mathcal{O}(1)$
and the piecewise linear Lagrange FEM might not converge.
However, if $k=2$, then
$\max_{K \in \T_h} (R_Kh_K) = \mathcal{O}(h^{2+\alpha - \beta})$,
and the finite element solution certainly converges to the exact
solution, although the convergence rate is worse than expected.
This means that ``bad'' triangulations with many very skinny triangles
can be remedied by using higher-order Lagrange elements.

\section{Numerical experiments}
\begin{figure}[thbp]
\begin{center}
      \includegraphics[width=5.2cm]{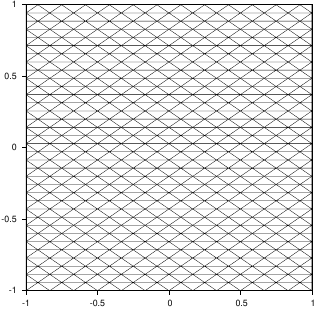}
 \qquad
          \includegraphics[width=7.2cm]{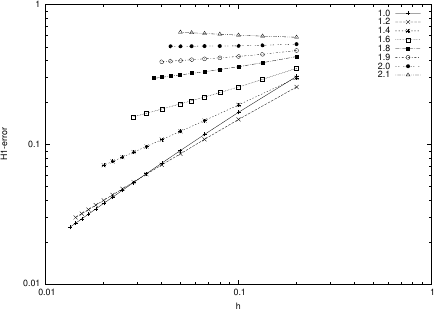}
\end{center}
 \caption{Triangulation of $\Omega$ with $N=12$ and $\alpha=1.6$,
 and the errors for FEM solutions in the $H^1$-norm.
 The horizontal axis represents the maximum diameter of the
 triangles and the vertical axis represents the $H^1$-norm of the errors
 of the FEM solutions. The number next to the symbol
indicates the value of $\alpha$.}
\label{fig3}
\end{figure}

To confirm the results obtained, we perform numerical experiments
similar to those in \cite{HKK}.  Let $\Omega := (-1,1)\times(-1,1)$, 
$f(x,y) := a^2/(a^2 - x^2)^{3/2}$, and
$g(x,y) := (a^2 - x^2)^{1/2}$ with $a:=1.1$.
Then we consider the following
Poisson equation:  Find $u \in H^1(\Omega)$ such that
\begin{equation}
   - \Delta u = f \quad \text{ in } \Omega, \qquad
   u = g \quad \text{ on } \partial\Omega.
  \label{test_problem}
\end{equation}
The exact solution of \eqref{test_problem} is $u(x,y) = g(x,y)$
and its graph is a part of the cylinder.
For a given positive integer $N$ and $\alpha > 1$, we consider the
isosceles triangle with base length $h:=2/N$ and height
$2/\lfloor 2/h^\alpha\rfloor \approx h^\alpha$, as shown
in Figure~\ref{fig3}.  Let $R$ be the circumradius of the
triangle.  For comparison, we also consider the isosceles
triangle with base length $h$ and height $h/2$ for 
$\alpha=1$.  We triangulate $\Omega$
with this triangle, as shown in Figure~\ref{fig3}. Let $\tau_h$
be the triangulation.  As usual, the set $S_h$ of piecewise linear
functions on $\tau_h$ and its subsets are defined by
\begin{align}
   S_h & := \left\{v_h \in C(\overline{\Omega}) 
  \bigm| v|_{K} \in \PP_1(K), \;
  \forall K \in \tau_h\right\}, \\
  S_{hg} & := \left\{v_h \in S_h 
  \bigm| v_h = g \; \text{at boundary nodes} \right\}, \\
   S_{h0} & := \left\{v_h \in S_h \mid v_h = 0 \text{ on }
 \partial \Omega\right\}.
\end{align}
Then, the piecewise linear finite element method for
 \eqref{test_problem} is defined as follows: Find $u_h \in S_{hg}$
such that
\begin{align*}
  (\nabla u_h, \nabla v_h)_\Omega = (f, v_h)_\Omega, \quad
  \forall v_h \in S_{h0},
\end{align*}
where $(\cdot, \cdot)_{\Omega}$ is the inner product of $L^2(\Omega)$.
By C\'ea's lemma and the result obtained, we obtain the estimation
\begin{align}
  |u - u_h|_{1,2,\Omega} \le
 \inf_{v_h \in S_{hg}} |u - v_h|_{1,2,\Omega} 
 \le \left(\sum_{K \in \tau_h}
   |u - \I_K^1 u|_{1,2,K}^2\right)^{1/2}
  \le CR|u|_{2,2,\Omega}.
\end{align}

The behavior of the error is given in Figure~\ref{fig3}.  The
horizontal axis represents the mesh size measured by the maximum
diameter of triangles in the meshes and the vertical axis represents the
error associated with FEM solutions in the $H^1$ semi-norm.  The graph
clearly shows that the convergence rates worsen as $\alpha$ approaches
$2.0$. For $\alpha=2.1$, the FEM solutions even diverge.  This is a
counterexample to the vaguely believed dogma that ``FEM solutions
always converge to the exact solution if $h \to 0$''. 
See also \cite{Oswald}.

\begin{figure}[thbp]
\begin{center}
    \includegraphics[width=7.2cm]{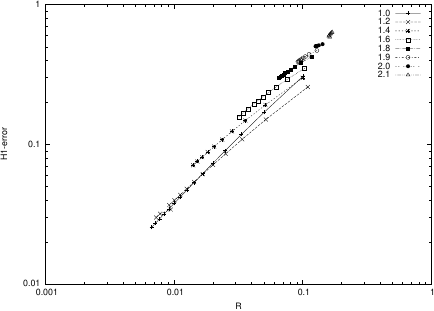}
\end{center}
 \caption{Replotted data: the errors in the $H^1$-norm of FEM
solutions measured using the circumradius. The horizontal axis represents
the maximum circumradius of the triangles.}
\label{fig4}
\end{figure}

We replot the
same data in Figure~\ref{fig4}, in which the horizontal axis represents the
maximum of the circumradius of triangles in the meshes.
Figure~\ref{fig4} shows convergence rates are almost the same in all
cases if we measure these with the circumradius.
These experiments strongly support that our theoretical results
are correct and optimal.

\vspace{5mm}

\newpage
\noindent
{\large \textbf{Acknowledgments}} \\
We thank Dr.\ Th\'eophile Chaumont-Frelet for his
valuable comments.

\end{document}